\documentclass{elsarticle}

\usepackage{mathtools}
\usepackage{amssymb}
\usepackage{siunitx}
\usepackage{bm}
\usepackage{bbm}
\usepackage{xcolor}
\usepackage{graphicx}
\usepackage{booktabs}
\usepackage[T1]{fontenc}
\usepackage[utf8]{inputenc}
\usepackage{subcaption}
\usepackage[margin=1in]{geometry}
\usepackage{tabularx}
\usepackage{booktabs}
\usepackage[colorlinks = true,
            linkcolor = blue,
            urlcolor  = blue,
            citecolor = blue,
            anchorcolor = blue,
            pdfauthor=author]{hyperref}
\usepackage[nameinlink]{cleveref}
\graphicspath{{../figures/pushcart/}{../figures/holzapfel_benchmark/}
  {../figures/holzapfel_benchmark/boxplots/}
  {../figures/holzapfel_benchmark/speedup/}
  {../figures/holzapfel_benchmark/boxplots/}
  {../figures/holzapfel_benchmark/data/middle_inner_1/}
  {../figures/holzapfel_benchmark/data/middle_outer_2/}
  {../figures/land_benchmark/}
  {../figures/land_benchmark/plots/figures/}
  {../figures/pushcart/stresses/}
  {../figures/pushcart/speedup/}
  {../figures/pushcart/violin/}
  {../figures/pushcart/}}

\newcommand{\tensor}[1]{\bm{#1}}             
\newcommand{\isotensor}[1]{\overline{\tensor{#1}}}             
\renewcommand{\vec}[1]{{\bm{#1}}}  
\newcommand{\Rvec}[1]{\underline{#1}}                
\newcommand{\Xvec}[1]{\bm{#1}} 
\newcommand{\Grad}{\operatorname{Grad}}      

\newcommand{\dx}{\,\mathrm{d}\vec{x}}
\newcommand{\dX}{\,\mathrm{d}\Xvec{X}}
\newcommand{\dsx}{\,\mathrm{d}s_{\vec{x}}}
\newcommand{\dsX}{\,\mathrm{d}s_{\Xvec{X}}}

\newcommand{\symotimes}{\overset{\mathrm{S}}{\otimes}}
\newcommand{\bb}[1]{\mathbbm{#1}}

\begin{document}
\begin{frontmatter}


\title{An accurate, robust, and efficient finite element framework for
       anisotropic, nearly and fully incompressible elasticity}
\author[1]{Elias Karabelas}
\author[2]{Matthias A. F. Gsell}
\author[2,3]{Gernot Plank}
\author[1,3]{Gundolf Haase}
\author[2]{Christoph M. Augustin\corref{cor1}}
\address[1]{%
  Institute for Mathematics and Scientific Computing,
  Karl-Franzens-University Graz, Graz, Austria}
\address[2]{%
  Gottfried Schatz Research Center: Division of Biophysics,
  Medical University of Graz, Graz, Austria}
\address[3]{%
  BioTechMed-Graz, Graz, Austria}
  \cortext[cor1]{Address correspondence to Christoph M. Augustin,
    Gottfried Schatz Research Center: Division of Biophysics,
    Medical University of Graz, Neue Stiftingtalstrasse 6/IV, Graz 8010, Austria.
  Email: christoph.augustin@medunigraz.at}

\begin{abstract}
  Fiber-reinforced soft biological tissues are typically
  modeled as hyperelastic, anisotropic, and nearly incompressible materials.
  To enforce incompressibility a multiplicative split of the deformation gradient
  into a volumetric and an isochoric part is a very common approach.
  However, due to the high stiffness of anisotropic materials in the preferred
  directions, the finite element analysis of such problems often suffers from
  severe locking effects and numerical instabilities.
  In this paper, we present novel methods to overcome locking phenomena for
  anisotropic materials using stabilized P1-P1 elements.
  We introduce different stabilization techniques and demonstrate the high
  robustness and computational efficiency of the chosen methods.
  In several benchmark problems we compare the approach to standard linear
  elements and show the accuracy and versatility of the methods to
  simulate anisotropic, nearly and fully incompressible materials.
  We are convinced that this numerical framework offers the possibility to
  accelerate accurate simulations of biological tissues, enabling
  patient-specfic parameterization studies, which require numerous forward simulations.
 \end{abstract}
 \journal{}

\end{frontmatter}
%
\section{Introduction}
Computer models of biological tissues, e.g., the simulation of vessel mechanics
or cardiac electro-mechanics (EM), aid in understanding the biomechanical
function of the organs and show promise to be a powerful tool for predicting
therapeutic responses. Advanced applications include the simulation of growth and remodeling
processes occurring in the failing heart or arteries~\cite{Grytsan2017,Genet2016,Peirlinck2019,Niestrawska2020}
as well as rupture risk assessment in
arterial aneurysms~\cite{gasser2010biomechanical,gultekin2016phase}.
Here, predictions of \emph{in-silico} models are often based on the
computation of local stresses, hence, an accurate computation of strain and
stress is indispensable to build confidence in simulation outcomes.
Additionally, high computational efficiency and excellent strong scaling
properties are crucial to perform simulations with highly-resolved, complex,
or heterogeneous geometries; to facilitate model personalization using a
large number of forward simulations; and to simulate tissue behavior over a
broad range of experimental protocols and extended observation periods.

\emph{In-silico} models of cardiac tissue and vessel walls
are typically based on the theory of hyperelasticity and properties of soft tissues
include a nonlinear relationship between stress and strain with large deformations
and a nearly-incompressible, anisotropic materials~\cite{Fung1990,guccione1991passive,holzapfel2000AS}.
Commonly, the resulting non-linear formulations are approximately solved using a
finite element (FE) approach~\cite{Augustin2014,Guccione1995,nash2000computational,sack2018construction}.
However, volumetric locking phenomena that are resulting in ill-conditioned
global stiffness matrices are frequently encountered.
In fact, this one of the
classical problems of modeling nearly incompressible hyperelasticity
\cite{babuska1992locking, hughes1987finite, zienkiewicz2000finite}.
Locking, often completely invalidating FE solutions, is in particular
prevalent for fiber-reinforced soft biological tissues due to a high
stiffness in the preferred fiber directions and thus
extensivley studied in recent publications
\cite{Farrell2021, gultekin2018quasi,hurtado2021accelerating, schroder2016}.

Typically, the modeling of (nearly) incompressible elastic materials
involves a split of the deformation gradient into a volumetric and an
isochoric part~\cite{flory1961thermodynamic}. Here, locking phenomena
are very common when using unstable approximation pairs such as
Q1-P0 or P1-P0 elements, i.e., when linear shape functions are the choice to
approximate the displacement field $\vec{u}$ and the hydrostatic pressure $p$
is statically condensed from the system of equations on the element-level.
It is well known that in such cases solution algorithms are likely to show very low
convergence rates, and that variables of interest such as stresses can be
inaccurate~\cite{gultekin2018quasi}.

To some degree locking problems in anisotropic hyperelasticity
for these simple elements can be reduced by using
augmented Lagrangian methods~\cite{glowinski1989augmented,simo1991quasi},
formulations with an unsplit deformation gradient for the anisotropic contribution
\cite{sansour2008physical, helfenstein2010nonphysical},
and methods with simplified kinematics for the anisotropic
contributions~\cite{schroeder2016novel}.
Another possibility to obtain more accurate results is the use of
higher order polynomials to approximate the
displacement~\cite{gerach2021electro,kerckhoffs2007coupling,
    kerckhoffs2008cardiac, Usyk2002}.
However, the incompressibility constraint is still modeled by a penalty
formulation, hence, volumetric locking may still be an issue and the modeling
of fully incompressible materials is not possible.
Additionally, already for quadratic ansatz functions the considerable
larger amount of degrees of freedom increases computational cost significantly.

A more sophisticated approach -- also allowing the modeling of fully incompressible materials --
is the reformulation of the underlying equations into a saddle point problem
by introducing the hydrostatic pressure $p$ as an additional unknown to the
system.
Here, from mathematical theory, approximation pairs for $\vec{u}$ and $p$ have
to fulfill the Ladyzhenskaya--Babu\^ska--Brezzi (LBB) or \emph{inf-sup} conditions
\cite{babuska1973finite,brezzi1974existence, chapelle1993infsup} to guarantee stability.
A popular choice are quadratic ansatz functions for the displacement and
linear ansatz functions for the pressure, i.e., the Taylor--Hood
element~\cite{taylor1973numerical,Nobile2012}.
Though stable, this element leads to a vast increase in degrees of freedom and
consequently a high computational burden; especially for applied problems in
the field of tissue mechanics with highly resolved geometries.

A computationally more favorable choice are equal order pairs with a
stabilization, widely used for linear and isotropic
elasticity~\cite{franca1988new, hughes1986new,masud2013framework, rossi2016implicite, chiumenti2015mixed, codina2000stabilization, lafontaine2015explicit}.
Yet, their extension to non-linear, anisotropic problems is
challenging~\cite{auricchio2005stability,auricchio2010importance,schroeder2016novel}.
In the specific case of modeling biological tissues Hu–Washizu-based
formulations are often used,
e.g.,~\cite{boerboom2003finite,goktepe2011computational,weiss1996finite,Zdunek2014}.
However, especially for problems undergoing large strains, this mixed three field
approach shows limited performance and robustness~\cite{schroeder2016novel}.

As our results suggest, a very promising and efficient stabilization approach for nearly
incompressible, anisotropic elasticity problems is a variant of the
MINI element~\cite{arnold1984},
originally established for computational fluid dynamics problems.
This element is modified for the application of incompressible hyperelasticity
and a bubble function is included in the space of displacements.
To improve efficiency, the support of this bubble can be eliminated from
the system of equation using static condensation.
First uses of MINI elements have been reported~\cite{ong2015,lee2017} though still using a piecewise constant ansatz for the hydrostatic pressure.
Even more efficient and notably simple to implement is a pressure projection
method originally introduced for the Stokes problem~\cite{dohrmann2004stabilized}.
To the best of our knowledge the here proposed methods were not
yet used in this form for anisotropic and nearly incompressible materials.

A big advantage of both stabilization techniques is that they do not rely on
artificial stabilization parameters that may influence the numerical solution.
We illustrate in different benchmarks that the same setting can be used
for a large variety of tissue mechanics problems allowing for a one-for-all approach.
By comparing to literature we show that our methods are suitable to
compute accurate strain and stress fields and outperform existing contributions
in terms of efficiency.

The paper is outlined as follows:
in \Cref{sec:methods} we recall the mathematical background of modeling
anisotropic, nearly incompressible elasticity and introduce the
theoretical framework of our stabilization techniques.
Subsequently, \Cref{sec:results} documents three benchmarks problems to show
the applicability of the stabilized P1-P1 elements in different scenarios.
For each benchmark we give a detailed problem description and discuss results
and computational efficiency by comparing to the literature and analyzing
strong-scaling properties.

To show the usefulness of the presented methods to clinically relevant problems,
we present a 3D EM model of the heart that is coupled to a 0D model of blood flow.
This constitutes the most complete model of cardiac EM in the literature to date
as all components, i.e., electrophysiology, cellular dynamics, active stress, passive tissue
mechanics, pre- and afterload, are based on physiological, state-of-the-art
models from the literature.
Finally, \Cref{sec:conclusion} concludes the paper with a brief summary and
all required equations to implement the methods in a software framework are
given in the Appendix.

\section{Methodology}\label{sec:methods}
\subsection{Almost Incompressible Nonlinear Elasticity}
Let \(\Omega_0 \subset \bb{R}^3\) denote the reference configuration and
let \(\Omega_t \subset \bb{R}^3\) denote the current configuration
of the domain of interest.
Assume that the boundary of \(\Omega_0\) is decomposed into
\(\partial \Omega_0 = \Gamma_{\mathrm{D},0} \cup \Gamma_{\mathrm{N},0}\)
with \(|\Gamma_{D,0}| > 0\).
Here, $\Gamma_{\mathrm{D},0}$ describes the Dirichlet part of the boundary
and $\Gamma_{\mathrm{N},0}$ describes the Neumann part of the boundary, respectively.
Further, let \(\vec{N}\) be the unit outward normal on \(\partial \Omega_0\).
The nonlinear mapping
\(\vec{\phi}\colon \Xvec{X} \in \Omega_0 \rightarrow \vec{x} \in \Omega_t\),
defined by
\(\vec{\phi}:= \Xvec{X} + \vec{u}(\Xvec{X},t)\), with displacement \(\vec{u}\),
maps points in the reference configuration to points in the current configuration.
Following standard notation, we introduce the
\emph{deformation gradient} \(\tensor{F}\), the Jacobian $J$,
and the \emph{left Cauchy-Green tensor} $\tensor{C}$ as
\begin{align*}
  \tensor{F} := \Grad \vec{\phi} = \tensor{I} + \Grad \vec{u},\quad
  J := \det(\tensor{F}),\quad
  \tensor{C} := \tensor{F}^\top \tensor{F}.
\end{align*}
Here, \(\Grad(\bullet)\) denotes the gradient with respect to the reference
coordinates \(\vec{X} \in \Omega_0\). The displacement field \(\vec{u}\)
is sought as infimizer of the functional
\begin{align}
  \nonumber\Pi^\mathrm{tot}(\vec{u}) &:= \Pi(\vec{u}) + \Pi^\mathrm{ext}(\vec{u}),\\
  \nonumber\Pi(\vec{u}) &:= \int\limits_{\Omega_0} \Psi(\tensor{F}(\vec{u}))\dX,\\
 \label{eq:potential_ext}
 \Pi^{\mathrm{ext}}(\vec{u})
 &:= - \rho_0\int\limits_{\Omega_0} \vec{f}(\vec{X})  \cdot \vec{u}\dX
 - \int\limits_{\Gamma_\mathrm{N,0}} \vec{h}(\vec{u}) \cdot \vec{u}\dsX,
\end{align}
over all admissible fields $\vec{u}$ with $\vec{u} = \vec{g}_\mathrm{D}$ on
$\Gamma_\mathrm{D,0}$, where, \(\Psi\) denotes the strain energy function;
\(\rho_0\) denotes the material density in reference configuration;
\(\vec{f}\) denotes a volumetric body force;
\(\vec{g}_\mathrm{D}\) denotes a given boundary displacement;
and \(\vec{h}\) denotes a given follower surface traction defined as
\begin{align*}
  \vec{h}(\vec{u}) := -p_\mathrm{ext} J(\vec{u}) \tensor{F}^{-\top}(\vec{u}) \vec{N},
\end{align*}
with giving external load $p_\mathrm{ext} > 0$.
For ease of presentation it is assumed that \(\rho_0\) is constant and
\(\vec{f}\), and \(\vec{g}_\mathrm{D}\) do not depend on \(\vec{u}\).
Existence of infimizers is, under suitable assumptions,
guaranteed by the pioneering works of Ball, see~\cite{Ball1976} and~\cite{ciarlet2002} for the case of follower loads.

In this study, we consider nearly incompressible materials,
meaning that \(J \approx 1\). A possibility to model this behavior was
originally proposed by~\cite{flory1961thermodynamic}
using a split of the deformation gradient \(\tensor{F}\) such that
\begin{equation}
  \tensor{F} = \tensor{F}_\mathrm{vol} \isotensor{F}.
\end{equation}
Here, \(\tensor{F}_\mathrm{vol}\) describes the volumetric change
while \(\isotensor{F}\) describes the isochoric change.
By setting \(\tensor{F}_\mathrm{vol} := J^{\frac{1}{3}} \tensor{I}\) and
\(\isotensor{F} := J^{-\frac{1}{3}} \tensor{F}\)
we get \(\det(\isotensor{F}) = 1\) and \(\det(\tensor{F}_\mathrm{vol}) = J\).
Analogously, by setting \(\tensor{C}_\mathrm{vol} := J^{\frac{2}{3}} \tensor{I}\) and
\(\isotensor{C} := J^{-\frac{2}{3}} \tensor{C}\),
we have \(\tensor{C} = \tensor{C}_\mathrm{vol} \isotensor{C}\).
Assuming a hyperelastic material, the Flory split also postulates an additive
decomposition of the strain energy function
\begin{equation} \label{eq:flory_strain_energy}
  \Psi = \Psi(\tensor{C}) =  U(J) + \overline{\Psi}(\isotensor{C}).
\end{equation}
The function \(U(J)\) will be used in the form
\begin{align*}
  U(J) := \frac{\kappa}{2} \Theta(J)^2
\end{align*}
where \(\kappa\) denotes the \emph{bulk modulus}.
In the literature many different choices for the functions \(\Theta(J)\) are proposed, see e.g \cite{rossi2016,hartmann2003polyconvexity,doll2000}
for examples and related discussion.
For studying also the limit case $\kappa \to \infty$ we will consider a reformulation
of~\Cref{eq:potential_ext} as perturbed Lagrangian-multiplier functional, see \cite{sussman1987finite,atluri1989formulation,simo1991,brink1996} for details.
Introducing the hydrostatic pressure $p$ we seek infimizers $(\vec{u}, p) \in V_{\vec{g}_{\mathrm{D}}} \times Q$ of
\begin{align}\label{eq:pertb_lagrange}
\Pi^\mathrm{tot}(\vec{u}, p)&:= \Pi^\mathrm{PL}(\vec{u}, p) + \Pi^\mathrm{ext}(\vec{u}),\\
\nonumber\Pi^\mathrm{PL}(\vec{u}, p) &:= \int\limits_{\Omega_0}\overline{\Psi}(\overline{\tensor{C}}) + p \Theta(J(\vec{u}))-\frac{1}{2\kappa}p^2\dX.
\end{align}
To guarantee well-definedness, we assume that
\begin{align*}
  V_{\vec{g}_{\mathrm{D}}} &:= \left\{ \vec{v} \in [H^1(\Omega_0)]^3 : \vec{v}\lvert_{\Gamma_{D,0}} = \vec{g}_{\mathrm{D}} \right\},
\end{align*}
with $H^1_0(\Omega_0)$ being the standard Sobolev space of square integrable functions
having a square integrable gradient, and $Q=L^2(\Omega_0)$.
For a more in-depth discussion we refer to \cite{Ball1976,ciarlet2002}.
To solve for infimizers of~\Cref{eq:pertb_lagrange} we calculate the variations with respect to $\vec{v}$ and $q$.
This results in the following non-linear variational problem, find $(\vec{u}, p) \in V_{\vec{g}_{\mathrm{D}}} \times Q$ such that
\begin{align}
  \label{eq:nonlinear_residual:1}R_\mathrm{vol}(\vec{u}, p; \vec{v}) &= 0, \\
  \label{eq:nonlinear_residual:2}R_\mathrm{inc}(\vec{u}, p; q) &= 0,
\end{align}
for all $(\vec{v},q) \in V_\mathrm{0} \times Q$.
Here,
\begin{align*}
  R_\mathrm{vol}(\vec{u}, p; \vec{v}) &:= a_\mathrm{isc}(\vec{u};  \vec{v}) + a_\mathrm{vol}(\vec{u}; \vec{v}) - l_\mathrm{follow}(\vec{u}, p_\mathrm{ext}; \vec{v}),\\
  R_\mathrm{inc}(\vec{u}, p; \vec{v}) &:= b_\mathrm{vol}(\vec{u};  q) - c(p, q),
\end{align*}
where
\begin{align*}
  a_\mathrm{isc}(\vec{u}; \vec{v}) &:= \int\limits_{\Omega_0} \tensor{S}_\mathrm{isc}(\vec{u}) : \tensor \Sigma(\vec{u}, \vec{v})\dX,
  &a_\mathrm{vol}(\vec{u}, p; \vec{v}) &:= \int\limits_{\Omega_0} p \tensor{S}_\mathrm{vol}(\vec{u}) : \tensor \Sigma(\vec{u}, \vec{v})\dX,\\
  b_\mathrm{vol}(\vec{u}; q) &:= \int\limits_{\Omega_0} \Theta(J(\vec{u})) q\dX,
  &c(p, q) &:= \frac{1}{\kappa} \int\limits_{\Omega_0} p q \dX,\\
  l_\mathrm{follow}(\vec{u}, p_\mathrm{ext}; \vec{v}) &:= -p_\mathrm{ext} \int\limits_{\Gamma_N} J(\vec{u}) \tensor{F}^{-\top} \Xvec{N} \cdot \vec{v}\dsX,
\end{align*}
with $\tensor \Sigma(\vec{u}, \vec{v}) := \mathrm{sym}(\tensor{F}^\top(\vec{u}) \Grad \vec{v})$.
Components of the second Piola--Kirchhoff stress tensor
\begin{equation}
  \tensor{S}_\mathrm{tot} = \tensor{S}_\mathrm{vol} + \tensor{S}_\mathrm{isc}
  \label{eq:stress_tensor}
\end{equation}
are computed as
\begin{align*}
 \tensor{S}_\mathrm{isc}
   :=  J^{-\frac{2}{3}} \mathrm{Dev}(\overline{\tensor{S}}), \quad
 \overline{\tensor{S}}
   := \frac{\partial\overline{\Psi}(\overline{\tensor{C}})}
           {\partial\overline{\tensor{C}}} \quad
 \tensor{S}_\mathrm{vol} := \pi(J) \tensor{C}^{-1}, \quad
 \pi(J) := J \Theta'(J).
\end{align*}
When modeling electrically active tissue, we consider an additive decomposition of the
isochoric part of the stress tensor.
The total stress tensor is now given by the additive decomposition
\begin{equation}\label{eq:total_stress}
\tensor{S}_\mathrm{tot} = \tensor{S}_{\mathrm{a}} + \tensor{S}_{\mathrm{p}} = \tensor{S}_{\mathrm{a}} + 2\frac{\partial{\rm \Psi}(\tensor{C})}{\partial \tensor{C}}.
\end{equation}

To simulate the effect of the circulatory system, these equations are coupled
to a 0D lumped model as in~\cite{augustin2021computationally}.
The corresponding nonlinear variational problem reads as find $(\vec{u}, p) \in V_{\vec{g}_{\mathrm{D}}} \times Q$ and $\Rvec{p}_{\mathrm{CAV}} \in \bb{R}^{n_\mathrm{CAV}}$ such that
\begin{align}
  \label{eq:nonlinear_residual_cv:1}R_\mathrm{vol}(\vec{u}, p, \Rvec{p}_{\mathrm{CAV}}; \vec{v}) &= 0, \\
  \label{eq:nonlinear_residual_cv:2}R_\mathrm{inc}(\vec{u}, p; q) &= 0,\\
  \label{eq:nonlinear_residual_cv:3}R_{\mathrm{CAV},i}(\vec{u}, p_{\mathrm{CAV},i}) &= 0,
\end{align}
for all $(\vec{v}, q) \in V_0 \times Q$, and $i=1,\ldots,n_\mathrm{CAV}$.
Here, the variations read as
\begin{align}
  R_\mathrm{vol}(\vec{u}, p, \Rvec{p}_\mathrm{CAV}; \vec{v}) &:= a_\mathrm{isc}(\vec{u}; \vec{v}) + a_\mathrm{vol}(\vec{u}; \vec{v})
  + \sum_{i=1}^{n_\mathrm{CAV}}l_\mathrm{follower}(\vec{u}, p_{\mathrm{CAV},i}; \vec{v}),\\
  R_\mathrm{inc}(\vec{u}, p; \vec{v}) &:= b_\mathrm{vol}(\vec{u}; q) - c(p, q),\\
  R_{\mathrm{CAV},i}(\vec{u}, p_{\mathrm{CAV},i}) &:= V_{\mathrm{CAV},i}(\vec{u}) - V_{\mathrm{CS}}(p_{\mathrm{CAV},i}),
\end{align}
where
\begin{align*}
  V_{\mathrm{CAV},i}(\vec{u}) &:= \frac{1}{3}\int\limits_{\Gamma_{\mathrm{CAV},i,0}} J (\Xvec{X} + \vec{u}) \cdot \tensor{F}^{-\top} \Xvec{N} \dsX,
\end{align*}
with $\Gamma_{\mathrm{CAV},i,0}$ denoting the closed surface of the $i^\text{th}$ cavity in reference configuration.
The expression for cavity volume $V_{\mathrm{CAV},i}$ follows from applying Nanson's formula to the definition of cavity volume in the current configuration
\begin{align*}
  V_{\mathrm{CAV},i} := \frac{1}{3}\int\limits_{\Gamma_{\mathrm{CAV},i}} \vec{x} \cdot \vec{n} \dsx.
\end{align*}
For a more detailed account on the coupling of nonlinear elastic equations with 0D lumped parameter models we refer to \cite{augustin2021computationally}.
\subsection{Consistent Linearization}
For the subsequent discretization we need the consistent linearization of~\eqref{eq:nonlinear_residual:1}--\eqref{eq:nonlinear_residual:2} and~\eqref{eq:nonlinear_residual_cv:1}--\eqref{eq:nonlinear_residual_cv:3}
and we obtain the following linear saddle-point problem: for each
\((\vec{u}^k, p^k) \in V_{\vec{g}_{\mathrm{D}}} \times Q\),
find \((\Delta \vec{u}, \Delta p) \in V_0 \times Q\) such that
\begin{align}
  \label{eq:saddle_point_nl:1}
  a_k(\Delta \vec{u}, \Delta \vec{v}) + a_{k,\Gamma}(\Delta \vec{u}, \Delta \vec{v}) + b_k(\Delta p, \Delta \vec{v})
  &= -R_\mathrm{vol}(\vec{u}^k, p^k; \Delta \vec{v}),\\
  \label{eq:saddle_point_nl:2}
  b_k(\Delta q, \Delta \vec{u}) - c(\Delta p, \Delta q)
  &= -R_\mathrm{inc}(\vec{u}^k, p^k;\Delta q),
\end{align}
where
\begin{align}
  \nonumber a_k(\Delta \vec{u}, \Delta \vec{v})
    &:= \int\limits_{\Omega_0} \Grad \Delta \vec{v} \tensor{S}_{\mathrm{tot},k}
    : \Grad \Delta \vec{u}\dX
    + \int\limits_{\Omega_0}\tensor \Sigma(\vec{u}_k, \Delta \vec{v})
    : \bb{C}_{\mathrm{tot},k} : \tensor \Sigma(\vec{u}_k, \Delta \vec{u})\dX,\\
  \label{eq:cons_lin_neumann}
  a_{k,\Gamma}(\Delta \vec{u}, \Delta \vec{v})
    &:= p_\mathrm{ext} \int\limits_{\Gamma_{N,0}} J_k (\tensor{F}^{-\top}_k
    : \Grad \Delta \vec{u}) \Delta \vec{v} \cdot \tensor{F}^{-\top} \Xvec{N} \dsX
   \\ \nonumber &- p_\mathrm{ext} \int\limits_{\Gamma_{N,0}} J_k \tensor{F}^{-\top}_k
    (\Grad \Delta \vec{u})^\top \tensor{F}^{-\top}_k \Xvec{N}
    \cdot \Delta \vec{v}\dsX,\\
  \nonumber b_k(\Delta p, \Delta \vec{v})
    &:= \int\limits_{\Omega_0} \Delta p \pi(J_k) \tensor{F}^{-\top}_k
    : \Grad \Delta \vec{v} \dX,
\end{align}
using the following abbreviations
\begin{align*}
  \tensor{F}_k &:= \tensor{F}(\vec{u}_k), &J_k &:= \det(\tensor{F}_k),\\
  \tensor{S}_{\mathrm{tot},k}
    &:= \left.\tensor{S}_{\mathrm{isc}}\right|_{\vec{u} =\vec{u}_k}
    + p_k \left. \tensor{S}_{\mathrm{vol}}\right|_{\vec{u} =\vec{u}_k},
  &\bb{C}_{\mathrm{tot},k}
    &:= \left.\bb{C}_{\mathrm{isc}}\right|_{\vec{u} =\vec{u}_k}
    + p_k \left.\bb{C}_{\mathrm{vol}}\right|_{\vec{u} = \vec{u}_k},\\
  \bb{C}_\mathrm{vol}
    &:= k(J) \tensor{C}^{-1} \otimes \tensor{C}^{-1}
    - 2 \pi(J) \tensor{C}^{-1} \odot \tensor{C}^{-1},
  &k(J) &:= J^2 \Theta''(J) + J \Theta'(J),
\end{align*}
and $\bb{C}_{\mathrm{isc}}$ given in \eqref{eq:def_c_isc}.
For the deviation of term \eqref{eq:cons_lin_neumann} see~\ref{sec:appendix_linearization},
other terms in \eqref{eq:saddle_point_nl:1}--\eqref{eq:saddle_point_nl:2} have been
discussed previously, see~\cite{karabelas2019MINI}.

In the case of an attached circulatory system we obtain the following linearized system, find \((\Delta \vec{u}, \Delta p, \Delta p_{\mathrm{CAV}})\) such that
\begin{align}
  a_k(\Delta \vec{u}, \Delta \vec{v}) + a_{k,\Gamma}(\Delta \vec{u}, \Delta \vec{v}) + b_k(\Delta p, \Delta \vec{v}) + l_{\mathrm{surface}}(u^k,\Delta p_{\mathrm{CAV}}; \Delta \vec{v}) &= -R_{\mathrm{vol}}(\vec{u}^k, p^k, p_{\mathrm{CAV}}^k; \Delta \vec{v}),\label{eq:linsys1}\\
  b_k(\Delta q, \Delta \vec{u}) - c(\Delta p, \Delta q) &= -R_{\mathrm{inc}}(\vec{u}^k, p^k; \Delta q), \label{eq:linsys2}\\
  d_k(\Delta \vec{u}) - e_k(\Delta p_{\mathrm{CAV}}) &= -R_{\mathrm{CAV}}(\vec{u}^k, p_\mathrm{CAV}^k), \label{eq:linsys3}
\end{align}
where
\begin{align}
  \label{eq:compl_matrix} e_k(\Delta p_{\mathrm{CAV}}) := \frac{\partial V_{\mathrm{CAV}}}{\partial p_{\mathrm{CAV}}},
\end{align}
and $d_k$ defined as in \eqref{eq:var_volume_cav}.
The term \eqref{eq:compl_matrix} depends on the chosen model for the circulatory system and a detailed discussion is out of the scope of this work.
For a detailed derivation of the explicit representation of the compliance matrix \eqref{eq:compl_matrix} stemming from the model used in Section \ref{sec:results} we refer to \cite{augustin2021computationally}.

\subsection{Finite Element Discretization}
Here we provide a summary of the finite element discretization used in the subsequent results.
The framework builds upon methods previously introduced for isotropic, passive
mechanics in~\cite{karabelas2019MINI}. In the following, we extend this approach to
anisotropic tissues also allowing for complex EM simulations that are
coupled to a 0D system of the circulatory system.

Let $\mathcal T_h$ be a finite element partitioning of $\overline{\Omega}$
consisting of tetrahedral and/or isoparametric hexahedral finite elements.
We assume standard regularity assumptions \cite{ciarlet2002} and invertibility of the isoparametric mapping $F_K$ from the reference element $\hat{K}$ to a physical element $K \in \mathcal T_h$.
For tetrahedral elements this poses no additional restrictions, for hexahedral elements we refer to \cite{knabner2003} for details.
Let further $\hat{P}_1$ and $\hat{Q}_1$ denote the space of lowest order linear/trilinear finite element functions on the reference tetrahedron/hexahedron.
The discrete analogue to~\eqref{eq:saddle_point_nl:1}--\eqref{eq:saddle_point_nl:2} reads as:
given $(\vec{u}_h^k, p_h^k) \in V_{h,\vec{g}_{\mathrm{D}}} \times Q_h$, find
$(\Delta \vec{u}_h, \Delta p_h) \in V_{h,0} \times Q_h$ such that
\begin{align}
  \label{eq:fem_saddle_point:1}
  a_k(\Delta \vec{u}_h, \vec{v}_h) + a_{k,\Gamma}(\Delta \vec{u}_h, \vec{v}_h) + b_k(\Delta p_h, \vec{v}_h) &= -R_\mathrm{vol}(\vec{u}_h^k,p_h^k; \vec{v}_h),\\
  \label{eq:fem_saddle_point:2}
  b_k(q_h, \Delta \vec{u}_h) - c(\Delta p_h, q_h) &= -R_\mathrm{inc}(\vec{u}_h^k, p_h^k; q_h)
\end{align}
for all $(\vec{v}_h, q_h) \in V_{h,0} \times Q_h$.
The discrete spaces $V_{h,0}$ and $Q_h$ are defined as
\begin{align*}
  V_{h,0} &:= \left\{\vec{v} \in \left[H^1_0(\Omega_0)\right]^3 : \vec{v} \vert_{K} = \hat{\vec{v}} \circ F_K^{-1}, \hat{\vec{v}} \in [\hat{\bb{V}}]^3,\forall K \in \mathcal T_h\right\},\\
  Q_h &:=\left\{q \in L^2(\Omega_0) : q \vert_{K} = \hat{q} \circ F_K^{-1}, \hat{q} \in \hat{\bb{Q}},\forall K \in \mathcal T_h\right\},
\end{align*}
and additionally we introduce
\begin{align*}
    V_{h,\vec{g}_{\mathrm{D}}} := H^1_{\vec{g}_{\mathrm{D}}}(\Omega_0) \cap V_{h,0}.
\end{align*}
\subsubsection{Pressure Projection Stabilized Equal Order Pair}
\label{sec:projection}
The pressure projection stabilization was originally introduced for solving
Stokes problems~\cite{dohrmann2004stabilized} and has also been applied in the context
of linear elasticity~\cite{cante2014,rodriguez2016}.
Recently, we extended its use to isotropic, nonlinear elasticity~\cite{karabelas2019MINI}.
A similar approach can be used for anisotropic materials, we set
$\hat{\bb{V}} := \hat{\bb{Q}} := \hat{P}_1/\hat{Q}_1$ for tetrahedral or hexahedral elements.
To ensure stability, we have to modify the definition of the residuals
in \eqref{eq:nonlinear_residual:2} and \eqref{eq:nonlinear_residual_cv:2} to
\begin{align}
  \nonumber
  \widetilde{R}_\mathrm{inc}(\vec{u}_h, p_h; q_h) &:=
    R_\mathrm{inc}(\vec{u}_h, p_h; q_h) - s_h(p_h, q_h), \\
  \label{eq:stab_def}
  s_h(p, q) &:=
    \int_{\Omega_0} \frac{1}{\mu_*} (p - \Pi_h p)(q - \Pi_h q)\dx,
\end{align}
where the projection operator $\Pi_h$ is defined elementwise as
\begin{align*}
    \Pi_h q \vert_{K} := \frac{1}{\lvert K \rvert }\int_{K} q\dx.
\end{align*}
In contrast to~\cite{karabelas2019MINI}, the parameter $\mu_*$ is no longer an
arbitrary value but set to $|K|^{1/3}$;
a choice that showed excellent results for all discussed
anisotropic problems in~\Cref{sec:results} as well as isotropic
benchmarks in~\cite{karabelas2019MINI}.
We note, that the integral in~\eqref{eq:stab_def} has to be understood as a sum over all the elements of the triangulation of the domain $\Omega_0$.
For a more comprehensive overview and implementation details we refer to \cite{karabelas2019MINI}.
\subsubsection{MINI Element}
\label{sec:mini_element}
One of the earliest strategies in constructing a stable finite element pairing for
discrete saddle-point problems arising from Stokes Equations is the MINI-Element, dating back to the works of
Brezzi et al, see for example \cite{arnold1984,brezzi1992}.
Briefly, the strategy is to enrich the basis of lowest order finite elements by adding a higher degree polynomial with support restricted to the interior of the element.
Thus, for the tetrahedral reference element $\hat{K}$ we define
\begin{align*}
  \hat{\bb{V}} &:= \hat{P}_1 \oplus \{\hat{\psi}_\mathrm{B}\} \\
  \hat{\bb{Q}} &:= \hat{P}_1,\\
  \hat{\psi}_\mathrm{B} &:= 256 \xi_0 \xi_1 \xi_2 (1-\xi_0 - \xi_1 - \xi_2),
\end{align*}
where $(\xi_0, \xi_1, \xi_2) \in \hat{K}$ see also \cite{boffi2013mixed}.

For the hexahedral reference element $\hat{K} = [-1,1]^3$ we define
\begin{align}
  \label{eq:def_hex_bubble}\hat{\bb{V}} &:= \hat{Q}_1 \oplus \{\hat{\psi}_{\mathrm{B},1}, \hat{\psi}_{\mathrm{B},2}\} \\
  \nonumber\hat{\bb{Q}} &:= \hat{Q}_1,\\
  \nonumber\hat{\psi}_{\mathrm{B},1} &:= (1-\xi_0)^2(1-\xi_1)^2(1-\xi_2)^2 \hat{\psi}_\alpha,\\
  \nonumber\hat{\psi}_{\mathrm{B},2} &:= (1-\xi_0)^2(1-\xi_1)^2(1-\xi_2)^2 \hat{\psi}_\beta,
\end{align}
for $(\xi_0, \xi_1, \xi_2) \in \hat{K}$ and $(\alpha,\beta) \in [1,8]$
denoting the indices of two ansatz functions for diagonal opposite nodes in $\hat{K}$,
see~\cite{karabelas2019MINI}.

Classical results \cite{boffi2013mixed}
guarantee the stability of the MINI-Element for tetrahedral meshes in the almost incompressible linear elastic case.
For hexahedral elements we were able to prove stability in the almost incompressible linear elastictiy case provided an enrichment like \eqref{eq:def_hex_bubble} of the displacement ansatz space by two bubble functions see \cite{karabelas2019MINI}.
Due to the compact support of the bubble functions, static condensation can be
applied to remove the interior degrees of freedom during assembly.
Static condensation can be done by standard procedures~\cite{boffi2013mixed}
with the exception of follower loads which is discussed in~\ref{sec:appendix_static_condensation}.
As a result, these degrees of freedom are not needed to be considered in the full global stiffness matrix
assembly which is a key advantage of the MINI element.

\subsection{Material Models}
\label{sec:materials}
Arterial and myocardial tissue as modeled in~\Cref{sec:results} is considered as a
non-linear, hyperelastic, nearly incompressible, and anisotropic material
with a layered organization of fibers.
To model this behaviour in our benchmark problems we used strain energy functions of the form~\eqref{eq:flory_strain_energy};
namely, the transverseley isotropic constitutive law
by~\citet{guccione1995finite}
\begin{equation}\label{eq:guccioneStrainEnergy}
  \Psi(\tensor{C}) =  U(J) + \overline{\Psi}(\isotensor{C})
  \quad \text{with}\quad
    \overline{\Psi}(\isotensor{C}) =
    \frac{a}{2}\left[\exp\mathcal{Q}(\isotensor{C}, \vec{a}_i)-1\right]
\end{equation}
with $a>\SI{0}{\kPa}$ a scaling parameter, $\vec{a}_i$ fiber directions,
and $\mathcal{Q}$ a function in terms of scalar strain components.
Further, we compared the standard formulation of a
separated Fung-type exponential model
\begin{equation}\label{eq:HolzapfelAS}
  \Psi(\tensor{C}) =  U(J) + \overline{\Psi}_\mathrm{AS}(\isotensor{C})
  \quad \text{with}\quad
    \overline{\Psi}_\mathrm{AS}(\isotensor{C}) =
     \overline{\Psi}_{\mathrm{iso}}(\isotensor{C})
    + \overline{\Psi}_{\mathrm{aniso}}(\isotensor{C}, \vec{a}_i)
\end{equation}
with the formulation using an unsplit deformation gradient for the
anisotropic contribution
\begin{equation}\label{eq:HolzapfelWAS}
  \Psi(\tensor{C}) =  U(J) + {\Psi}_\mathrm{WAS}(\tensor{C})\quad \text{with}\quad
    {\Psi}_\mathrm{WAS}(\tensor{C}) =
       \overline{\Psi}_{\mathrm{iso}}(\isotensor{C})
      + \Psi_{\mathrm{aniso}}(\tensor{C}, \vec{a}_i),
\end{equation}
introduced in \cite{sansour2008WAS, helfenstein2010WAS}.
The specific form of the volumetric, $U(J)$, isotropic,
$\overline{\Psi}_{\mathrm{iso}}$, and anisotropic,
$\overline{\Psi}_{\mathrm{aniso}}$/$\Psi_{\mathrm{aniso}}$,
contributions will be discussed later for each of the benchmark problems.
\section{Numerical Examples}
\label{sec:results}
%
Biomechanical applications often require highly resolved meshes and thus
efficient and massively parallel solution algorithms for the
linearized system of equations become an important factor to deal with the
resulting computational load.
Extending our previous implementations for cardiac EM~\cite{augustin2016anatomically} we
used the software \textit{Cardiac Arrhythmia Research Package} (CARP)
\cite{vigmond2008solvers} which makes use of the MPI based library \emph{PETSc}
\cite{petsc-user-ref}. We solve the stabilized saddle-point
problem~\eqref{eq:fem_saddle_point:1}--\eqref{eq:fem_saddle_point:2}
by using a GMRES method with a block preconditioner based on a
smoothed aggregation algebraic multigrid (GAMG) approach which is
incorporated in \emph{PETSc}~\cite{may2016extreme}.

In all of the following benchmark problems our goal was to study the performance
and accuracy of different finite element discretizations, namely
i) \emph{Q1/P1-P0-AS}: discretization with piecewise linear displacements and piecewise
constant pressure using the strain energy function~\eqref{eq:HolzapfelAS};
ii) \emph{Q1/P1-P0-WAS}: discretization with piecewise linear displacement and piecewise
constant pressure using the strain energy function~\eqref{eq:HolzapfelWAS};
iii) \emph{Projection}: equal order discretization with
piecewise linear displacements and pressure, stabilized as described
in~\Cref{sec:projection} using the strain energy function~\eqref{eq:HolzapfelAS};
iv) \emph{MINI}: discretization using MINI elements as described
in~\Cref{sec:mini_element} using the strain energy function~\eqref{eq:HolzapfelAS}.

\subsection{Extension, Inflation and Torsion of a Simplified Artery Model}%
\begin{figure*}[htbp]\centering
  \includegraphics[width=0.8\linewidth]{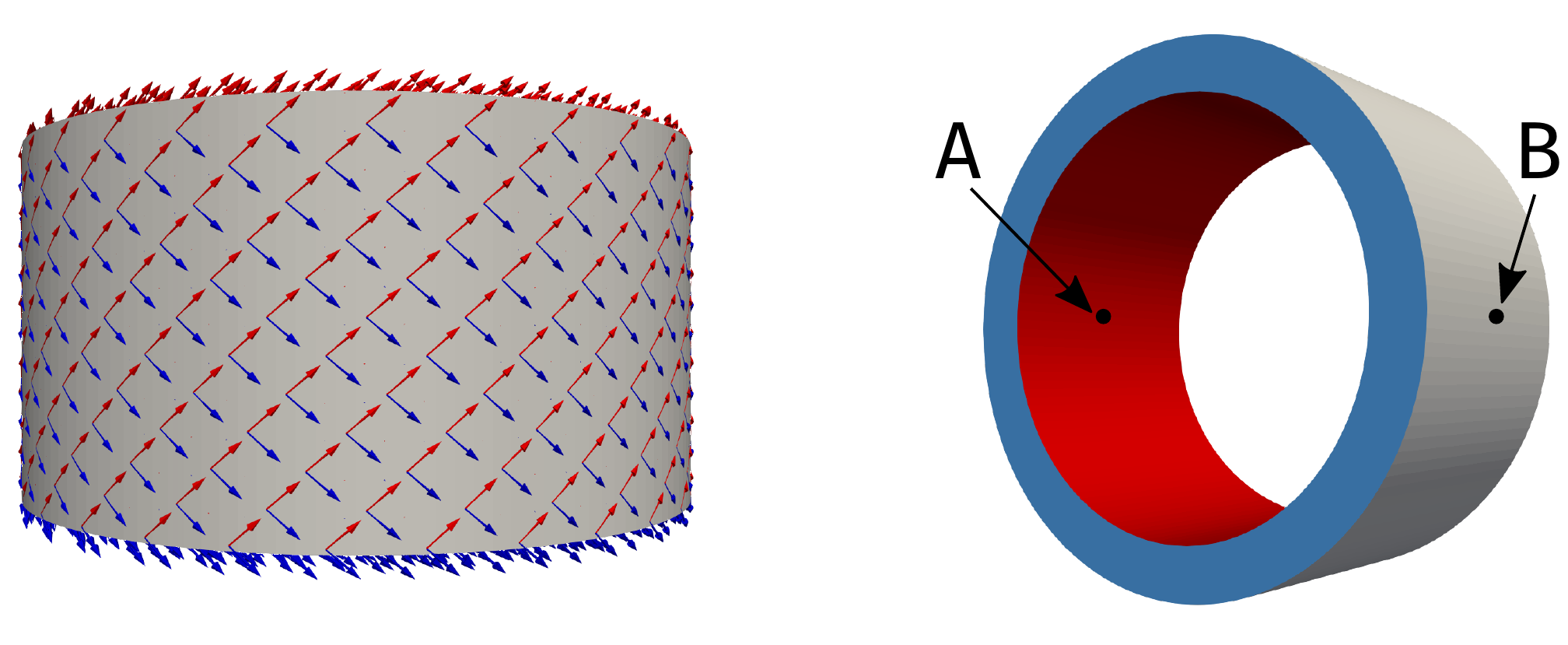}\\[-1em]
  {\small {(a)}}
  \hspace*{20em}
  {\small {(b)}}
  \caption{\emph{Artery benchmark:} Geometry, fiber distribution, and boundary conditions. Figure (a) shows the fiber arrangement and Figure (b) shows the geometry and colorcoded boundaries for the application of the boundary conditions and the evaluation points A, and B for mesh convergence.}%
  \label{fig:holzapfel_geo}
\end{figure*}

\label{sec:holzapfel_test}
\paragraph{Simulation setup}
First, we show the applicability of our proposed methods to a benchmark problem from~\citet{gueltekin2019}
where a simplified artery model is represented by a thick-walled cylindrical tube.
The dimensions of this idealized geometry with its centerline on the z-axis are as follows:
height $H=\SI{10}{\milli\meter}$, inner radius $R_1=\SI{8}{\milli\meter}$, and outer radius $R_2=\SI{10}{\milli\meter}$.
Two symmetric families of fibers, $\vec{f}_0$ and $\vec{s}_0$ are immersed in the tissue,
having an angle of $\SI{40}{\degree}$ with circumferential $\theta$-axis, see~\Cref{fig:holzapfel_geo}A.
As for loading, a monotonically increasing displacement up to $\SI{2}{\milli\meter}$
superimposed by a monotonically increasing torsion up to $\SI{60}{\degree}$
is applied on the top of the tube (marked blue in~\Cref{fig:holzapfel_geo}B).
Additionally, a linearly increasing pressure (follower load) up to $\SI{500}{\mmHg}$ is applied
on the inside of the tube (marked red in~\Cref{fig:holzapfel_geo}B).
Finally, the lower part of the tube is clamped at zero displacement.

The material is described by the strain-energy function~\eqref{eq:HolzapfelAS},
$\overline{\Psi}_\mathrm{AS}$, with
\begin{align*}
  \Theta(J) &:= J-1,\qquad
  \overline{\Psi}_{\mathrm{iso}}(\isotensor{C})
    := \frac{\mu}{2}\left(\overline{I}_1-3\right),\\
  \overline{\Psi}_{\mathrm{aniso}}(\isotensor{C},\vec{f}_0, \vec{s}_0)
  &:= \frac{k_1}{2 k_2}\sum_{i=4,6}
  \left(\mathrm{exp}\left(k_2\left(\overline{I}_i-1\right)^2\right)-1\right),
\end{align*}
and invariants
\[
\overline{I}_1 := \mathrm{tr}(\isotensor C),\quad
\overline{I}_4 := \isotensor C : \vec{f}_0 \otimes \vec{f}_0,\quad
\overline{I}_6 := \isotensor C : \vec{s}_0 \otimes \vec{s}_0,
\]
and analogously with $\Psi_{\mathrm{aniso}}(\tensor{C},\vec{f}_0, \vec{s}_0)$ for the $\mathrm{WAS}$ formulation,
$\Psi_\mathrm{WAS}$,~\eqref{eq:HolzapfelWAS}.
Material parameters were taken from~\cite[Table 1]{gueltekin2019}, i.e.,
$\kappa=\SI{5000}{\kPa}$, $\mu = \SI{10}{\kPa}$, $k_1=\SI{500}{\kPa}$, and $k_2=2.0$.
In case of the stabilized equal-order elements (Projection and MINI),
we set $1/\kappa = 0$ to render the material incompressible.
To assess mesh convergence simulations were performed on seven discretization levels, see~\Cref{tab:artery_meshes}.

\begin{figure*}[htbp]
  \centering
   \begin{subfigure}{0.47\textwidth}
  \includegraphics[width=\textwidth, keepaspectratio]{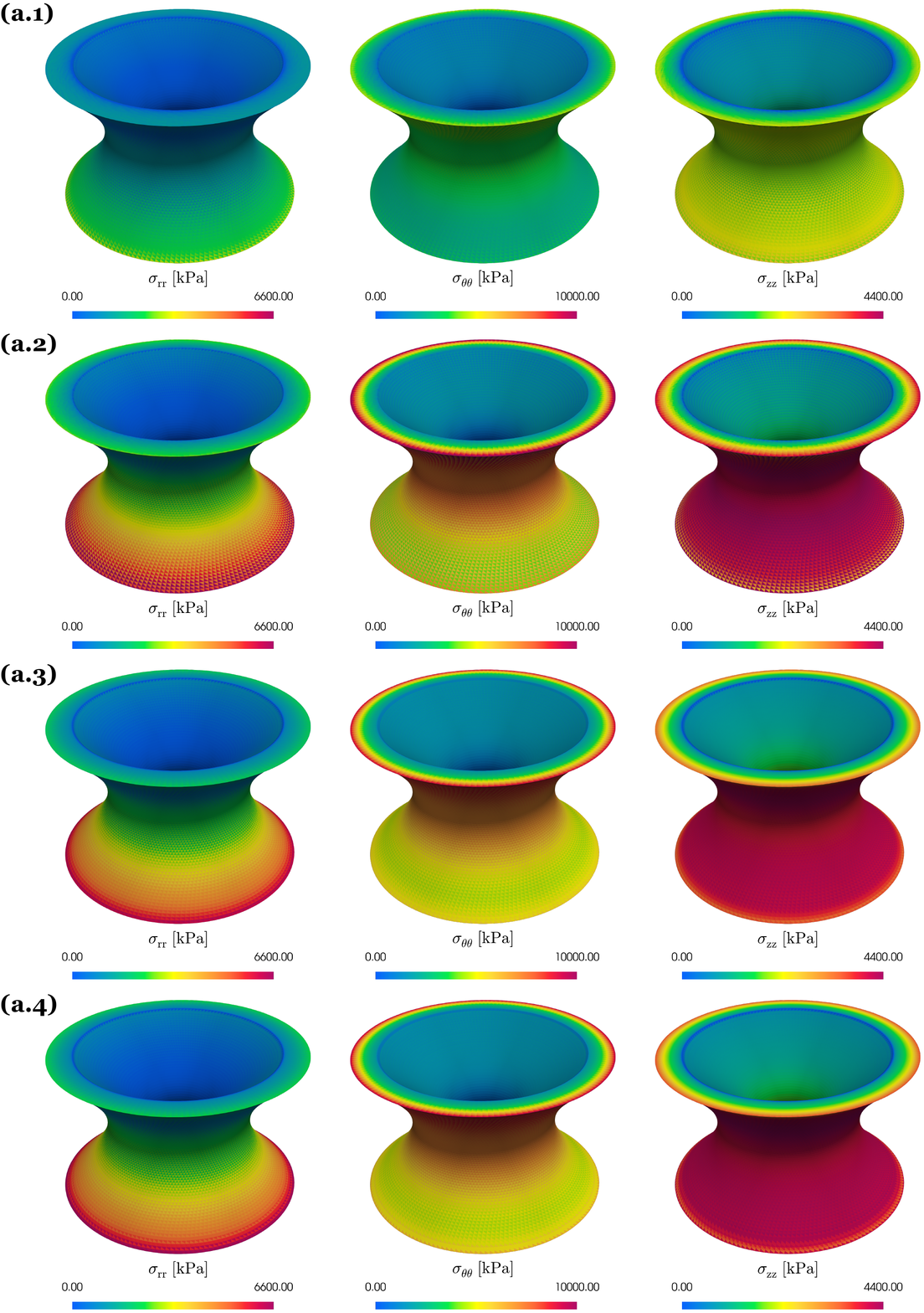}
  \subcaption{}
  \end{subfigure}\hspace{2em}
   \begin{subfigure}{0.47\textwidth}
  \includegraphics[width=\textwidth, keepaspectratio]{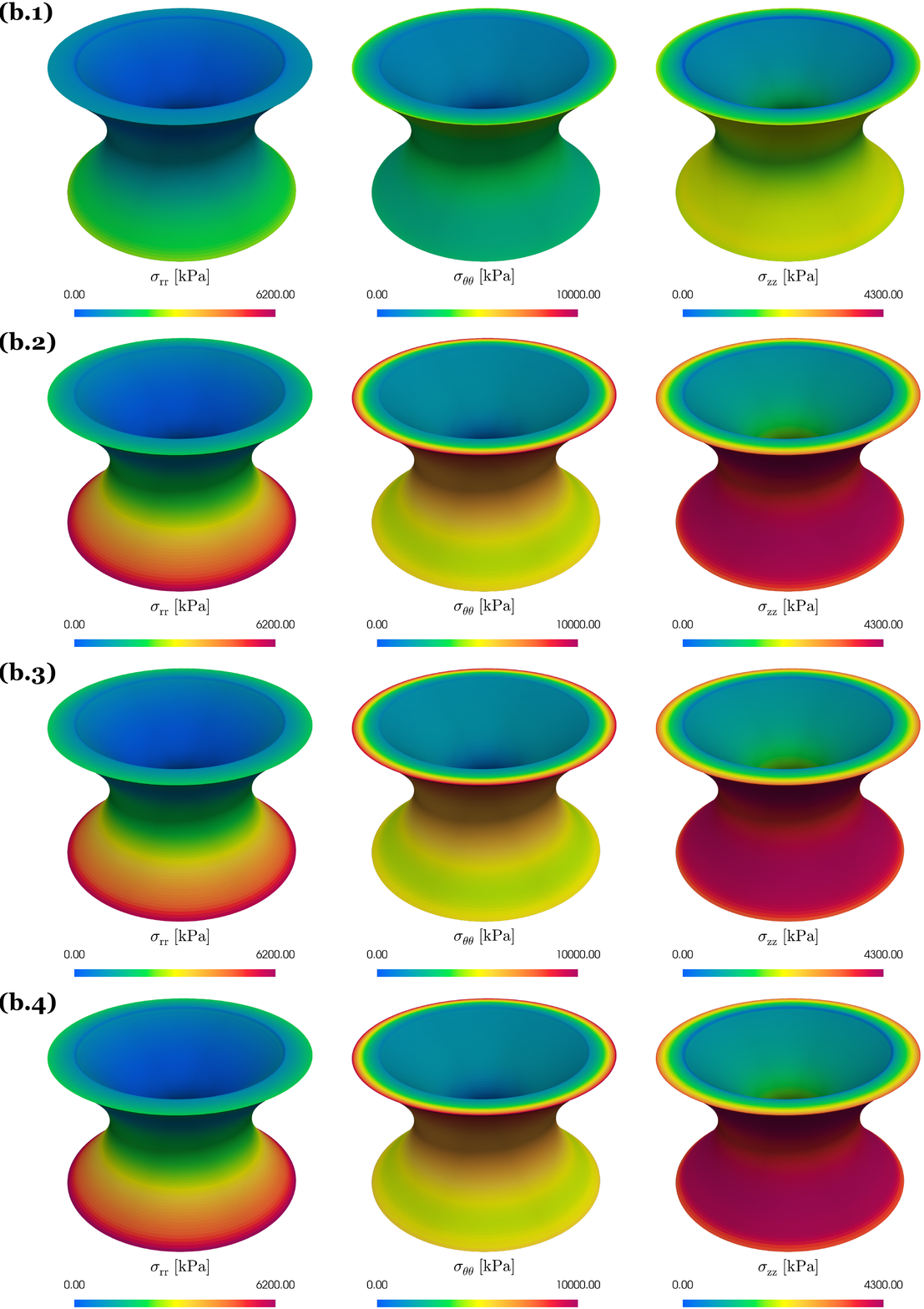}
  \subcaption{}
  \end{subfigure}
  \caption{\emph{Artery benchmark:} Stress distribution for (a) tetrahedral and
    (b) hexahedral elements.
    Shown are the the radial, $\sigma_\mathrm{rr}$, circumferential, $\sigma_{\theta \theta}$,
    and longitudinal stresses, $\sigma_{zz}$, on the
    finest discretization level $\ell=7$.
    The rows 1, 2, 3, 4 correspond to Q1-P0-AS, Q1-P0-WAS, Projection, and MINI,
    respectively.}%
  \label{fig:stress_pic}
\end{figure*}

\begin{figure*}[htbp]
  \centering
   \begin{subfigure}{0.85\textwidth}
  \includegraphics[width=1.\linewidth,keepaspectratio]{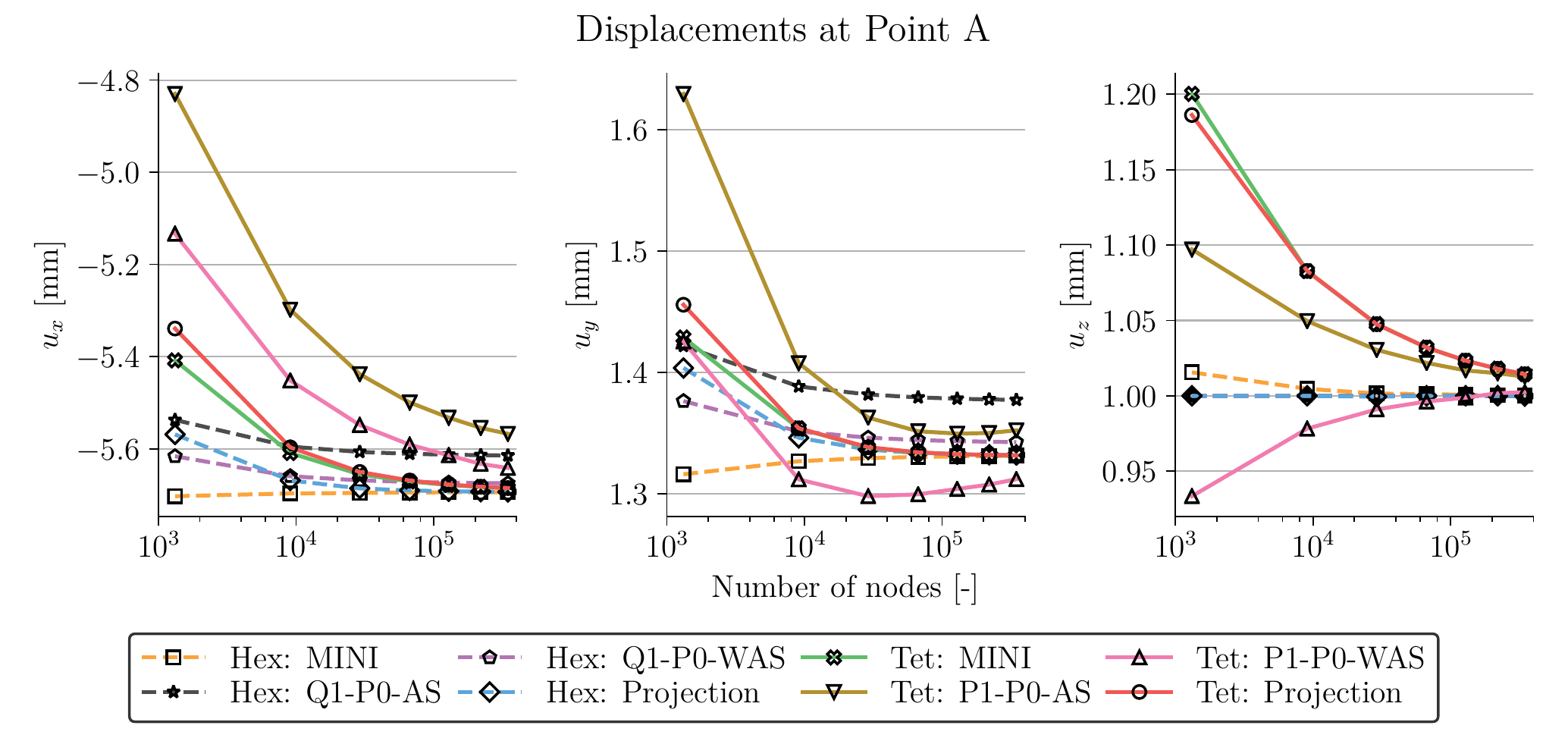}
  \subcaption{}
  \end{subfigure}
   \begin{subfigure}{0.85\textwidth}
  \includegraphics[width=1.\linewidth,keepaspectratio]{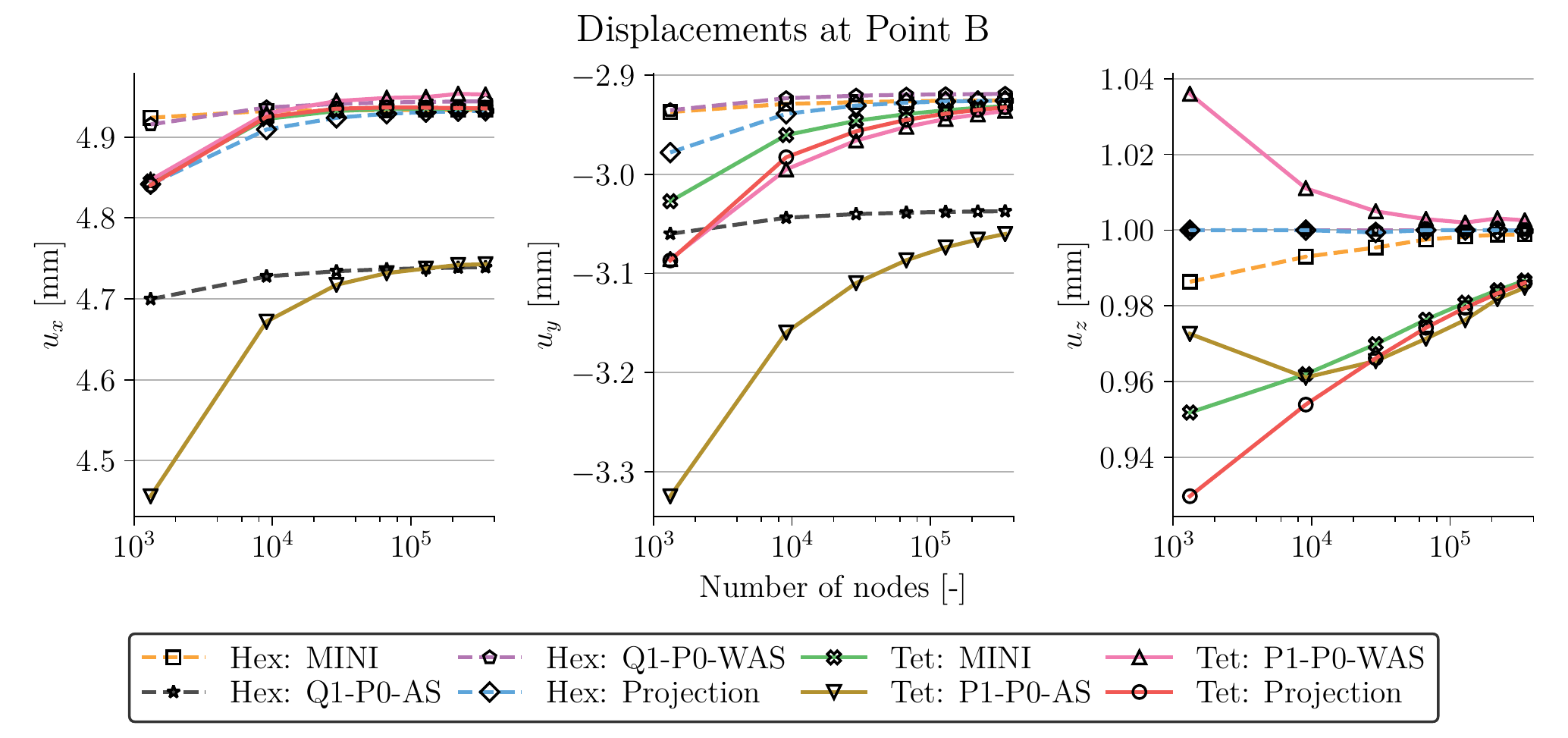}
  \subcaption{}
  \end{subfigure}
  \caption{\emph{Artery benchmark:} Mesh convergence. Shown are the individual displacement values at
    (a) Point A and (b) Point B for increasing mesh resolution and finite elements.}%
  \label{fig:disp_comp}
\end{figure*}

\begin{figure*}[htbp]
  \centering
   \begin{subfigure}{0.85\textwidth}
  \includegraphics[width=1.\linewidth,keepaspectratio]{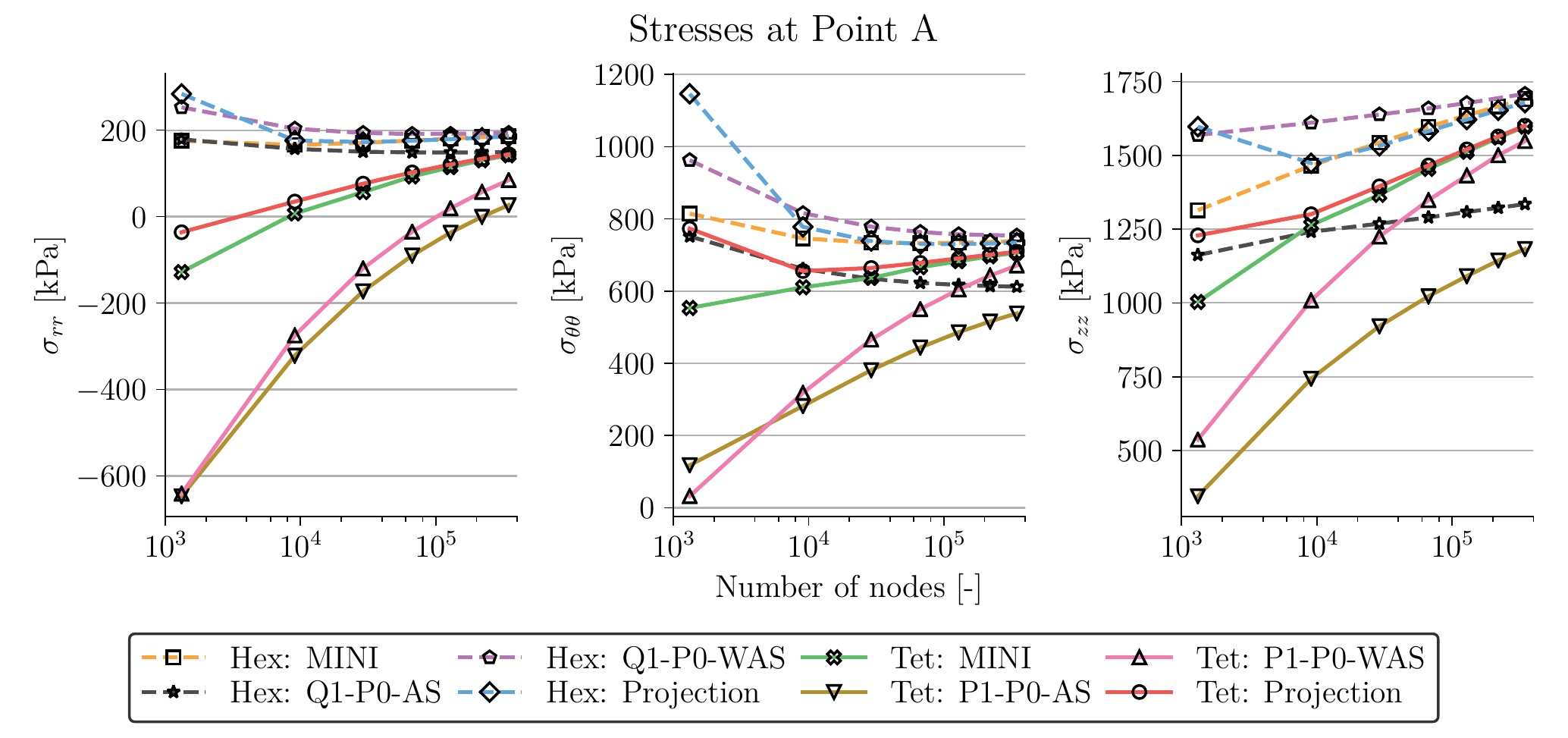}
  \subcaption{}
  \end{subfigure}
   \begin{subfigure}{0.85\textwidth}
  \includegraphics[width=1.\linewidth,keepaspectratio]{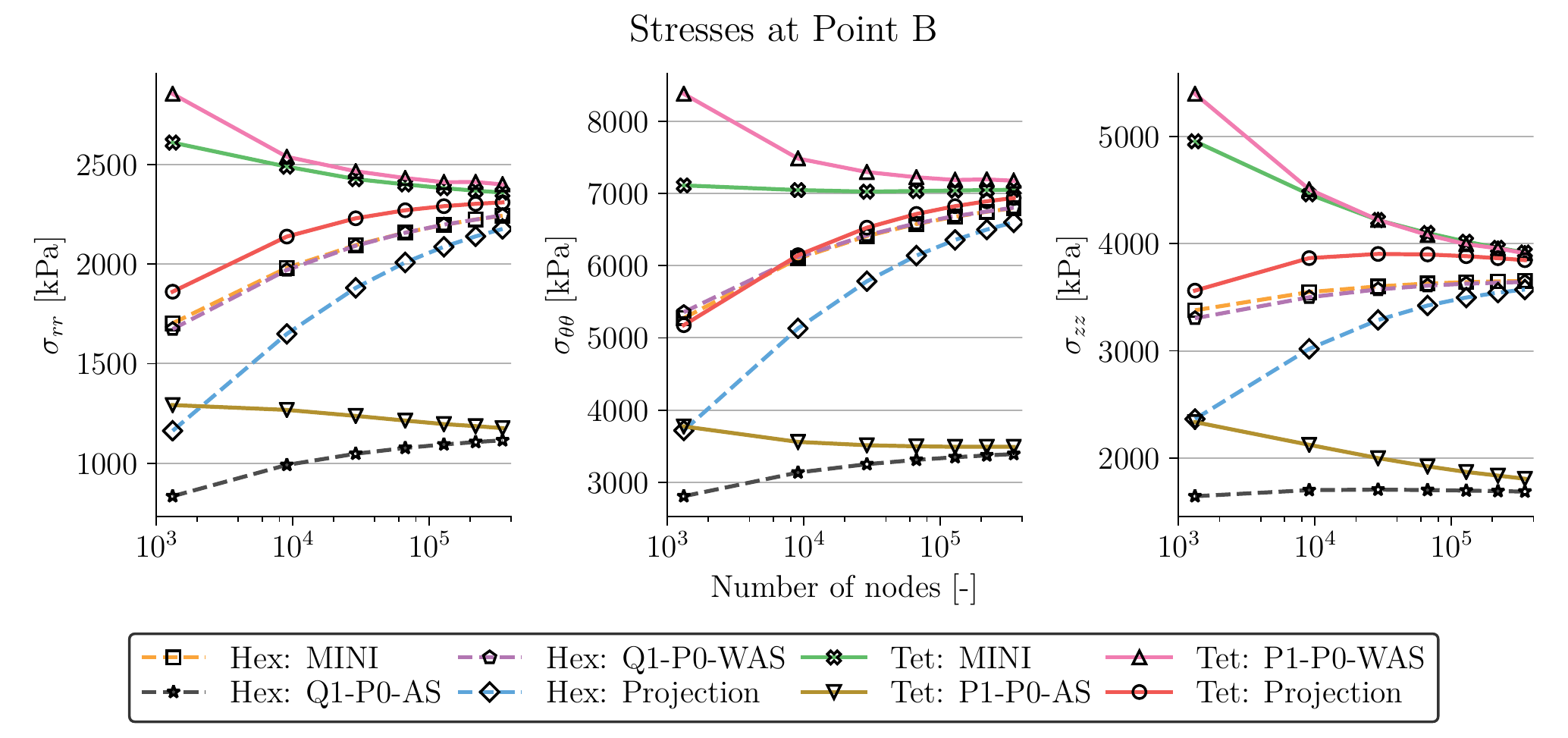}
  \subcaption{}
  \end{subfigure}
  \caption{\emph{Artery benchmark:} Mesh convergence. Shown are the individual stresses
    at (a) Point A and (b) Point B for increasing mesh resolution and finite elements.}%
  \label{fig:stress_comp}
\end{figure*}

\paragraph{Results} A comparison of the radial, $\sigma_\mathrm{rr}$, the circumferential,
$\sigma_{\theta \theta}$, and the axial, $\sigma_\mathrm{zz},$ components of the
Cauchy stress tensor is shown in~\Cref{fig:stress_pic} for the finest discretization level $\ell = 7$.
We see that with the exception of the lowest order discretizations with anisotropic splitting
(Q1/P1-P0-AS) the stress distribution is very similar and also matches results in~\citet{gueltekin2019}.
The observation that simulations with Q1/P1-P0-AS are not accurate for this benchmark problem is further
emphasized in~\Cref{fig:disp_comp,fig:stress_comp}.
Here,~\Cref{fig:disp_comp} shows the displacements $(u_x, u_y, u_z)$
and~\Cref{fig:stress_comp} shows the stress components
$(\sigma_\mathrm{rr}, \sigma_{\theta \theta}, \sigma_\mathrm{zz})$ at the evaluation points A and B over the discretization levels.
In agreement with~\cite{gueltekin2019} the lowest order discretization with
anisotropic splitting converges to a lower value than the other discretization types hinting at possible locking phenomena.
All other formulations perform similarly well. Here, discrepancies at the finest level $\ell=7$ are rather due to
differences in the meshes for tetrahedral and hexahedral grids.

\begin{figure*}[htbp]
  \centering
   \begin{subfigure}{0.75\textwidth}
  \includegraphics[width=1.\linewidth,keepaspectratio]{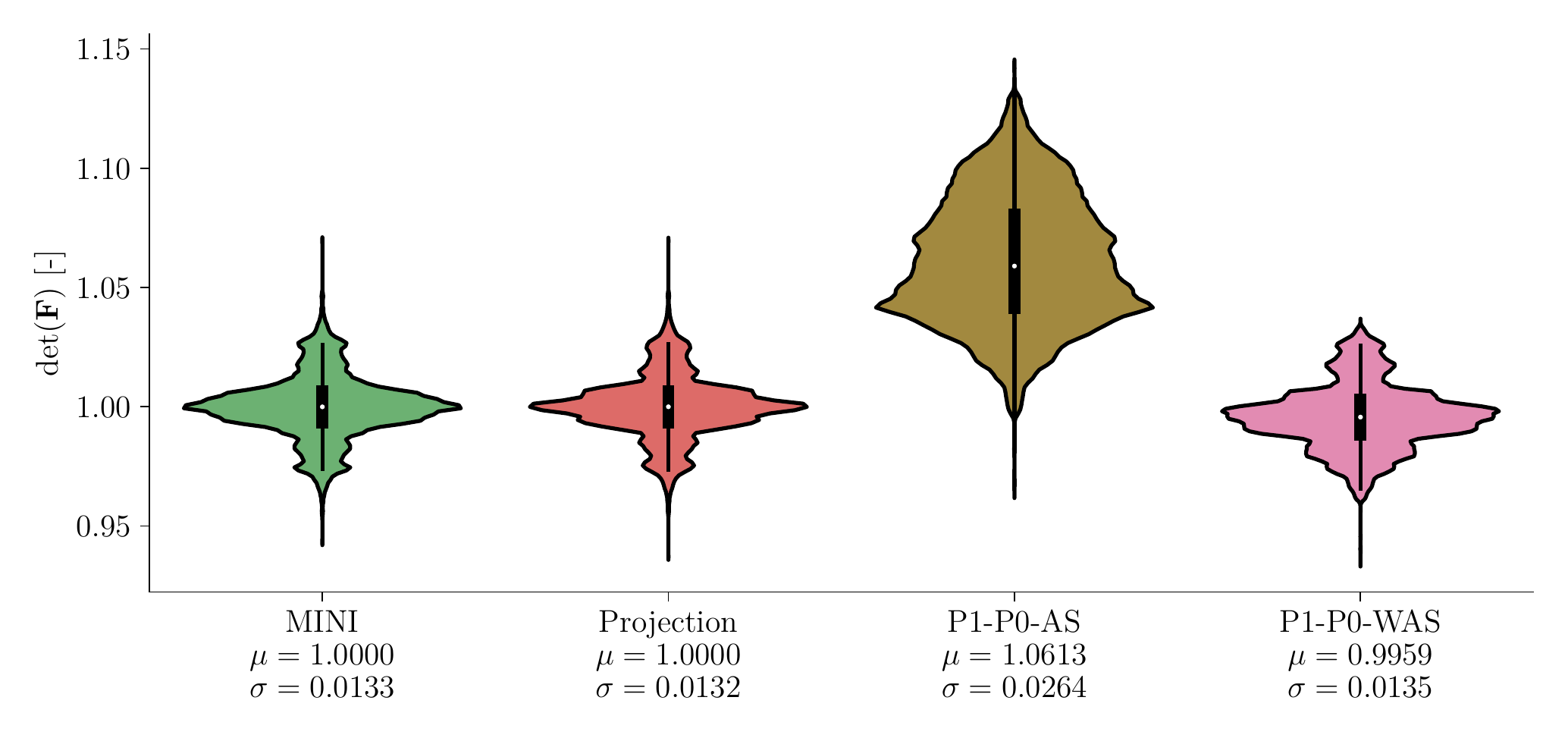}
  \subcaption{}
  \end{subfigure}
   \begin{subfigure}{0.75\textwidth}
  \includegraphics[width=1.\linewidth,keepaspectratio]{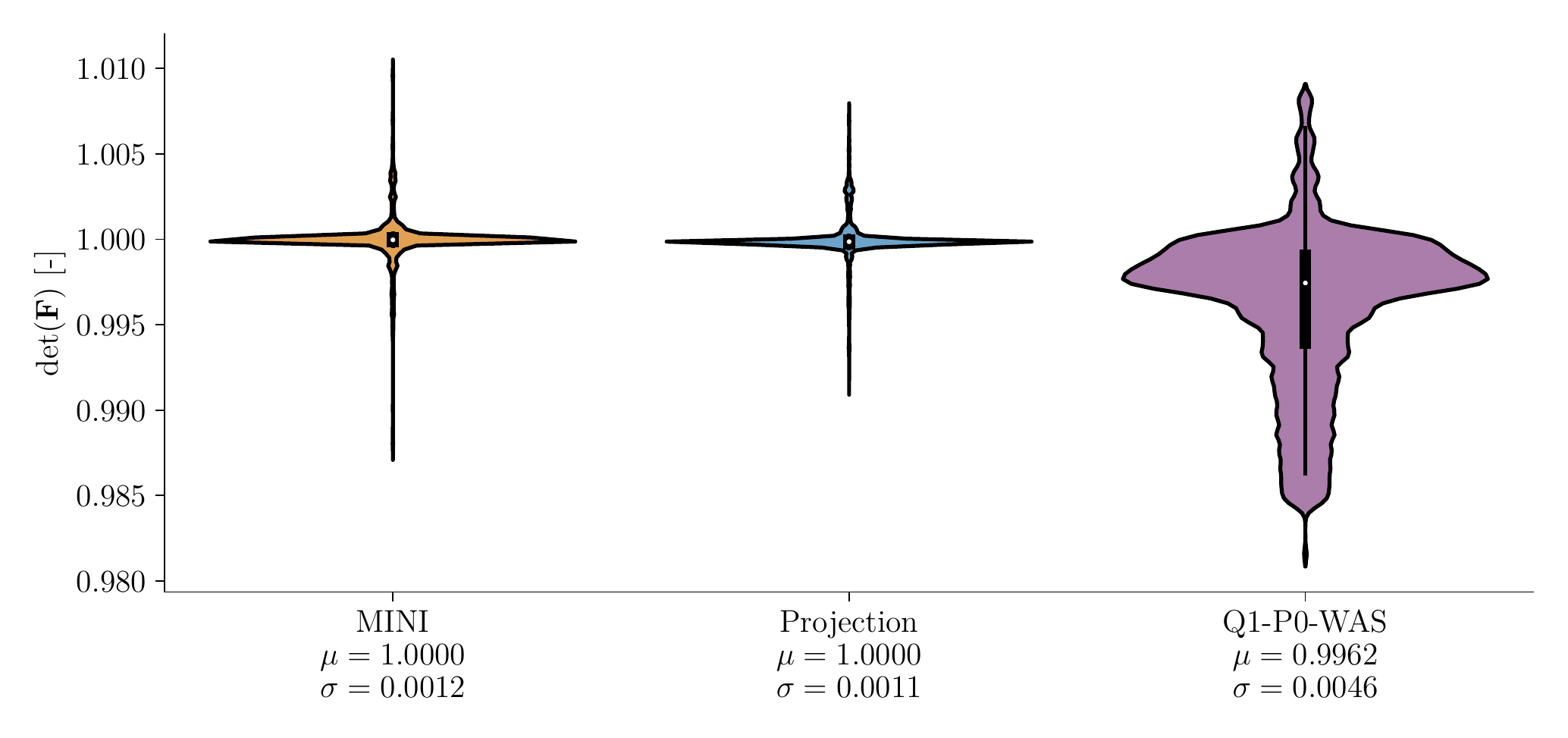}
  \subcaption{}
  \end{subfigure}
  \caption{\emph{Artery benchmark:} Jacobian distribution for (a) tetrahedral and
    (b) hexahedral elements.
    Shown are violin plots of the Jacobian distribution $\mathrm{det}(\tensor F)$
    on the finest discretization level at maximum loading for the various
    finite elements.
    Additionally, the mean $\mu$ and standard deviation $\sigma$ is given.
    Q1-P0-AS has been excluded from plot (b) as values were significantly higher
    compared to other element types.}%
  \label{fig:jacobian_comp}
\end{figure*}

\Cref{fig:jacobian_comp} shows a distribution of the Jacobian $\det(\tensor{F})$ on the finest level $\ell=7$.
Unsurprisingly, with a mean value $\mu$ close to 1 the saddle-point formulations (Projection, MINI) satisfy incompressibility
better than the penalty formulations (P1/Q1-P0-AS, P1/Q1-P0-WAS).
While the AS formulation led to a small increase in volume ($\mu>1$), the
WAS formulation resulted in a slightly reduced volume ($\mu<1$).

\begin{figure*}[htbp]
  \centering
  \includegraphics[width=0.85\linewidth,keepaspectratio]{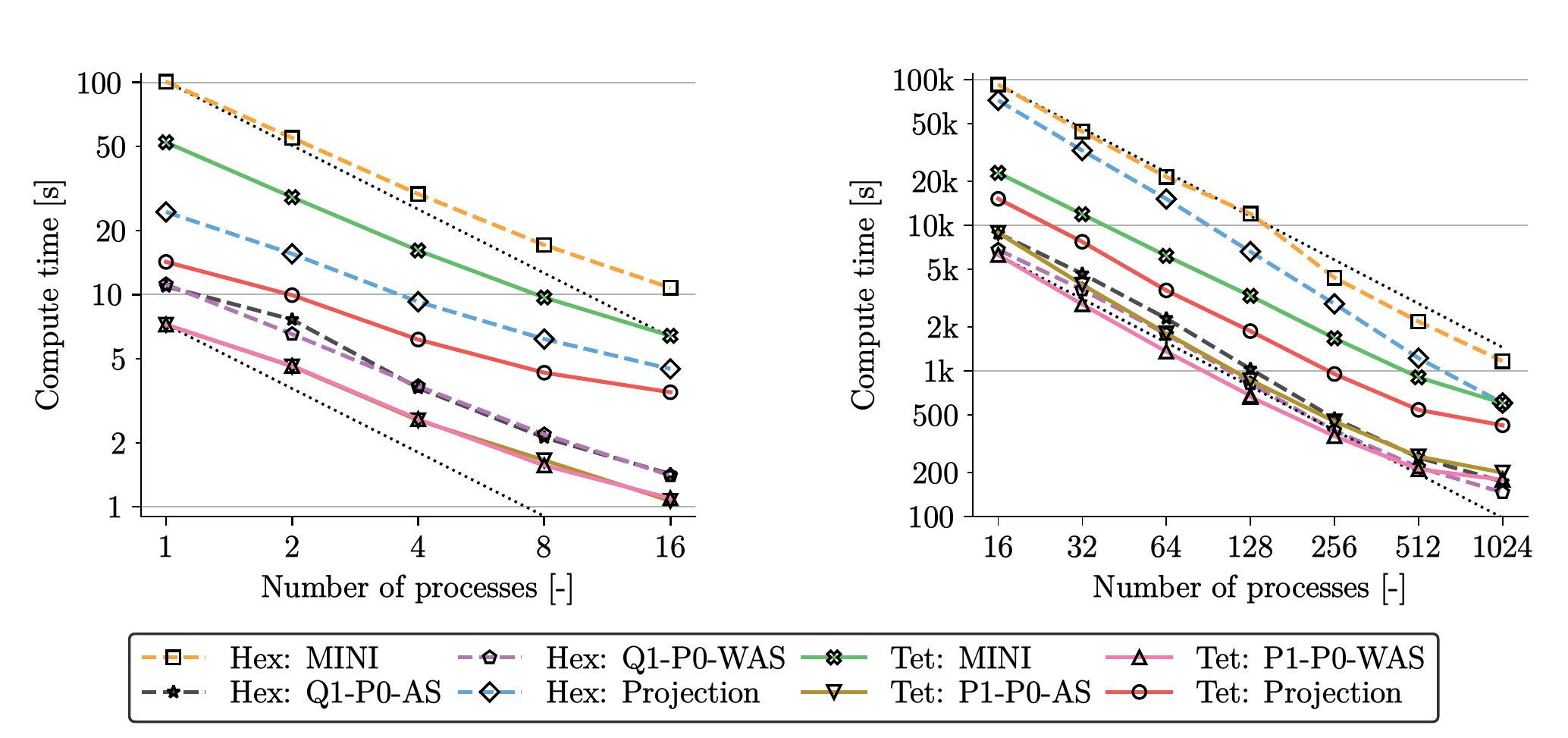}
  \caption{\emph{Artery benchmark:} Strong scaling results for the different element types for the coarsest grid ($\ell=1$, left)
    and the finest grid ($\ell=7$, right).
  Simulations were performed on 1 to 16 cores of a standard desktop computer for the coarse problem and
  on 16 to 1024 cores on Archer2.}%
  \label{fig:ab:strong_scaling}
\end{figure*}

\paragraph{Numerical performance} Computational times for the simulation using
different element types are given in~\Cref{fig:ab:strong_scaling};
left, for the  coarse problem ($\ell=1$) and right, for the finest grid ($\ell=7$).
For all cases we used a relative error reduction of $\epsilon=10^{-8}$ for the
GMRES linear solver and a relative error reduction of $\epsilon=10^{-6}$ for the residual of the Newton method.
Using a load stepping scheme, this required a total number of \num{100} linear solving steps for the coarse
problem and \num{1000} linear solving steps for the fine problem
to arrive at the final prescribed displacement and inner pressure of $\SI{500}{\mmHg}$.

Strong scaling was achieved for the coarse problem up to \num{16} cores on a desktop machine (AMD Ryzen Threadripper 2990X),
see~\Cref{fig:ab:strong_scaling}, left.
Here, computational times were averaged over \num{5} runs using the same setup and range from
\SIrange{1}{1.5}{\s} for penalty formulations (P1/Q1-P0-AS, P1/Q1-P0-WAS).
For the projection-stabilized element compute times were around \num{3} times slower
(Tet Projection: \SI{3.5}{\s}, Hex Projection: \SI{4.5}{\s}) while for
the MINI element compute times were around \num{7} times slower (Tet MINI: \SI{6.4}{\s}, Hex MINI: \SI{10.8}{\s}).

For the fine grid strong scaling was obtained up to \num{1024} cores on Archer2 (\url{https://www.archer2.ac.uk/}), see~\Cref{fig:ab:strong_scaling}, right.
Computational times using \num{1024} cores were between \SI{146}{\s} and \SI{200}{\s} for penalty formulations (P1/Q1-P0-AS, P1/Q1-P0-WAS);
around \num{3} times slower for the projection-stabilized elements (Tet Projection: \SI{422}{\s}, Hex Projection: \SI{601}{\s});
  and around \num{6} times slower for MINI elements (Tet MINI: \SI{600}{\s}, Hex MINI: \SI{1166}{\s}).

\begin{table}[htbp]
  \caption{\emph{Artery benchmark:} properties of idealized artery meshes used
           in~\Cref{sec:holzapfel_test}.}%
  \label{tab:artery_meshes}
  \centering
  \begin{tabular}{rrrr}
    \toprule
    $\ell$ & Elements (Hex) & Elements (Tet) & Nodes \\
    \midrule
    \num{1} & \num{960}    & \num{5760}    & \num{1320}  \\
    \num{2} & \num{7680}   & \num{46080}   & \num{9072} \\
    \num{3} & \num{25920}  & \num{155520}  & \num{29016} \\
    \num{4} & \num{61440}  & \num{368640}  & \num{66912} \\
    \num{5} & \num{120000} & \num{720000}  & \num{128520} \\
    \num{6} & \num{207360} & \num{1244160} & \num{219600} \\
    \num{7} & \num{329280} & \num{1975690} & \num{345912} \\
    \bottomrule
  \end{tabular}
\end{table}
\paragraph{Discussion}
In accordance with~\citet{gueltekin2019} we have shown for this benchmark that
the concept Q1/P1-P0-WAS is able to match quasi-incompressible responses compared to the gold standard of locking-free elements. For this extreme loading case
standard Q1/P1-P0 elements cannot reproduce accurate stress distributions;
not even on very fine grids, see~\Cref{fig:stress_comp}(b).

Further, this benchmark highlights the high efficiency of our methods.
While the authors in~\cite{gueltekin2019} reported computational times on
one 3.2GHz CPU unit that were 1 minute for Q1-P0-WAS elements and 20 minutes for their gold-standard (Q1-P0 elements using an Augmented Lagrangian method) we have
achieved compute times on a comparable 3.0GHz CPU unit and the same mesh that were
11 seconds for Q1-P0-WAS elements and 24 seconds for
locking-free Hex Projection elements, see~\Cref{fig:ab:strong_scaling}, left.
That means that stabilization techniques presented in this paper allow up to 50
times faster execution times for this benchmark compared to other gold standard methods.

Due to a higher number of linear iterations and higher matrix assembly times,
simulations with hexahedral meshes were more expensive compared
to simulations with tetrahedral grids. However, as also observed, e.g., by
\citet{chamberland2010comparison} hexahedral elements were slightly
more accurate than their tetrahedral equivalent.

\subsection{Inflation and Active Contraction of a Simplified Ventricle}%
\label{sec:ellipsoid}
\paragraph{Simulation setup}
To verify our EM setup we repeated the inflation and
active contraction benchmark from~\citet{land2015verification}.
We generated the reference geometry of an idealized ventricle
as a tetrahedral mesh of a truncated ellipsoid and prescribed a local orthonormal coordinate system with
fiber, $\vec{f}_0$, sheet, $\vec{s}_0$, and sheet-normal, $\vec{n}_0$, directions
according to this paper. We constructed three levels of refinement, see
see~\Cref{tab:ellipsoid_meshes} for discretization details.

\begin{table}[htbp]
  \caption{\emph{Ellipsoid benchmark:} properties of idealized
           \emph{ventricle meshes} used in~\Cref{sec:ellipsoid}.}%
  \label{tab:ellipsoid_meshes}
  \centering
  \begin{tabular}{rrr}
    \toprule
    $\ell$ & Elements & Nodes \\
    \midrule
    1 & \num{20709}    & \num{4732}  \\
    2 & \num{201495}   & \num{38896}  \\
    3 & \num{1747845}  & \num{312784} \\
    \bottomrule
  \end{tabular}
\end{table}

As material we used the transversely isotropic law
by~\citet{guccione1995finite}, see~\Cref{eq:guccioneStrainEnergy},
with $\Theta(J) := \mathrm{ln}(J)$ and
\begin{align*}
      \mathcal{Q} :=
        b_{\mathrm{f}} (\vec{f}_0\cdot\isotensor{E}\vec{f}_0)^2 +
        b_{\mathrm{t}} \left[(\vec{s}_0\cdot\isotensor{E}\vec{s}_0)^2+
                              (\vec{n}_0\cdot\isotensor{E}\vec{n}_0)^2+
                              2(\vec{s}_0\cdot\isotensor{E}\vec{n}_0)^2\right]+
        \quad 2b_{\mathrm{fs}} \left[(\vec{f}_0\cdot\isotensor{E}\vec{s}_0)^2+
        (\vec{f}_0\cdot\isotensor{E}\vec{n}_0)^2\right].
\end{align*}
Constitutive parameters were $a=\SI{2}{\kPa}$, $b_\mathrm{f} = 8$,
$b_\mathrm{t} = 2$, and $b_\mathrm{fs} = 4$.
In the benchmark paper the material is considered to be fully incompressible, hence,
 we chose $1/\kappa=0$ for the saddle-point formulation (Projection, MINI).
For the penalty formulation (P1-P0 elements) we chose $\kappa=\SI{1000}{\kPa}$
which was the best trade-off between convergence of the solver for all three levels
in~\Cref{tab:ellipsoid_meshes}, near incompressibility, and
minimization of locking effects.

As the material law above does not allow for a WAS formulation we repeated the
benchmark using a separated Fung-type exponential model
as in~\Cref{eq:HolzapfelAS}.
In particular, we chose a Holzapfel--Ogden material~\cite{Holzapfel2009} of the form
\begin{align}
  \Theta(J) &:= \ln(J),\qquad
  \overline{\Psi}_{\mathrm{iso}}(\isotensor{C})
    := \frac{a}{2b}\left\{\exp\left[b(\overline{I}_1-3)\right]-1\right\} \nonumber\\
  \overline{\Psi}_{\mathrm{aniso}}(\isotensor{C},\vec{f}_0, \vec{s}_0, \vec{n}_0)
   &:= \sum_{i=\mathrm{f,n}}\frac{a_i}{2b_i}\left\{\exp\left[b_i
      (\overline{I}_{4i}-1)^2\right]-1\right\}
    + \frac{a_\mathrm{fs}}{2b_\mathrm{fs}}\left\{\exp\left[b_\mathrm{fs}
    (\overline{I}_\mathrm{8fs})^2\right]-1\right\}, \label{eq:benchmark2_ortho}
\end{align}
with invariants
\[
\overline{I}_1 := \mathrm{tr}(\isotensor C),\quad
\overline{I}_{4\mathrm{f}}= \max\left(\vec{f}_0\cdot\isotensor{C}\vec{f}_0,1\right),\quad
\overline{I}_{4\mathrm{n}}= \max\left(\vec{n}_0\cdot\isotensor{C}\vec{n}_0,1\right),\quad
\]
such that contributions of compressed fibers are excluded,
and the interaction-invariant
\[
  \overline{I}_{8\mathrm{fs}}= \vec{f}_0\cdot\isotensor{C}\vec{s}_0.
\]
Analogously, we used the constitutive equation above with
$\Psi_{\mathrm{aniso}}(\tensor{C},\vec{f}_0, \vec{s}_0, \vec{n}_0)$ for the $\mathrm{WAS}$ formulation.
Material parameters were taken from~\cite{Guan2019AIC},
$a=\SI{0.809}{\kPa}$, $b=7.474$, $a_{\mathrm{f}}=\SI{1.911}{\kPa}$,
$b_{\mathrm{f}}=22.063$, $a_{\mathrm{n}}=\SI{0.227}{\kPa}$, $b_{\mathrm{n}}=34.802$,
$a_{\mathrm{fs}}=\SI{0.547}{\kPa}$, and $b_{\mathrm{fs}}=5.691$,
fitted to human myocardial experiments in~\cite{Sommer2015biomechanical}.

\begin{figure*}[htbp]
  \includegraphics[width=0.42\textwidth]{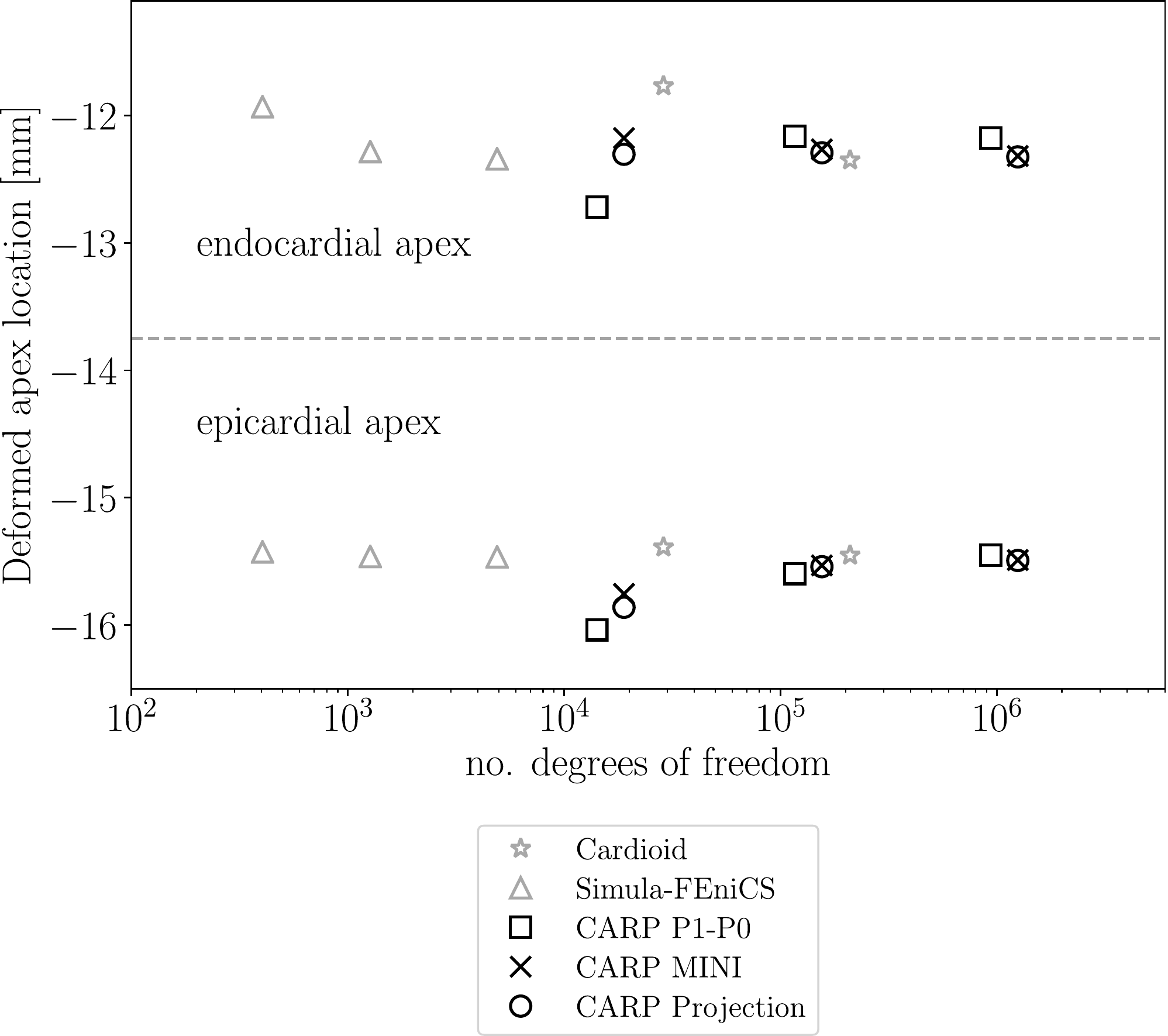}\hspace{4em}
  \includegraphics[width=0.42\textwidth]{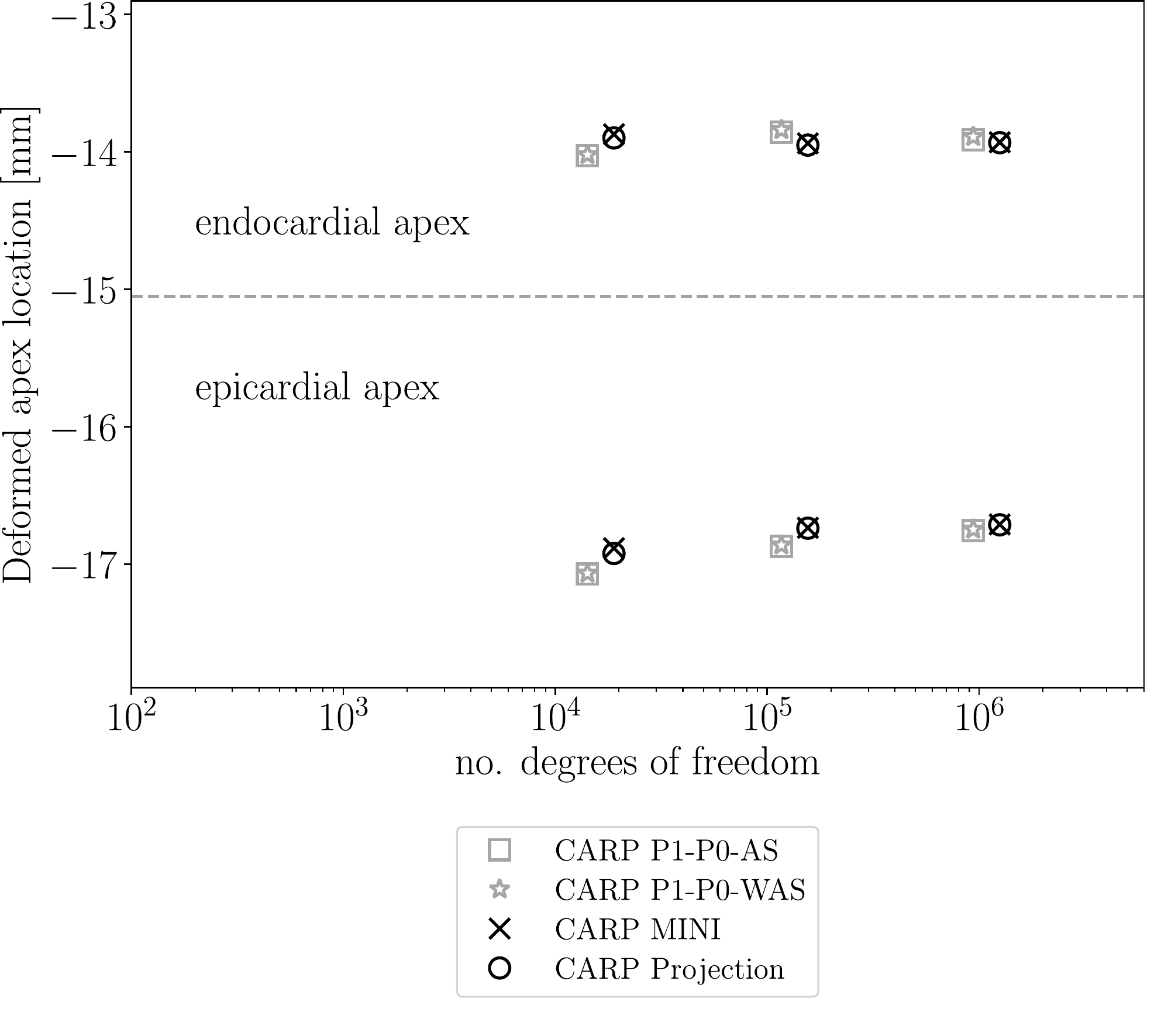}\\[-2.5em]
  \hspace*{1em}
  {\small {(a)}}
  \hspace*{21.8em}
  {\small {(b)}}\\[-1em]
  \caption{\emph{Ellipsoid benchmark:} apex location. The dashed line separates results for the deformed positions of the apex at the endo- and epicardium.
(a) Guccione material with comparison to benchmark results (in gray)
    presented in~\cite{land2015verification};
(b) Holzapfel--Ogden material with comparison to P1-P0-WAS formulation.}
  \label{fig:ellipsoid_apex}
\end{figure*}

\begin{figure*}[htbp]
  \includegraphics[width=1.0\linewidth]{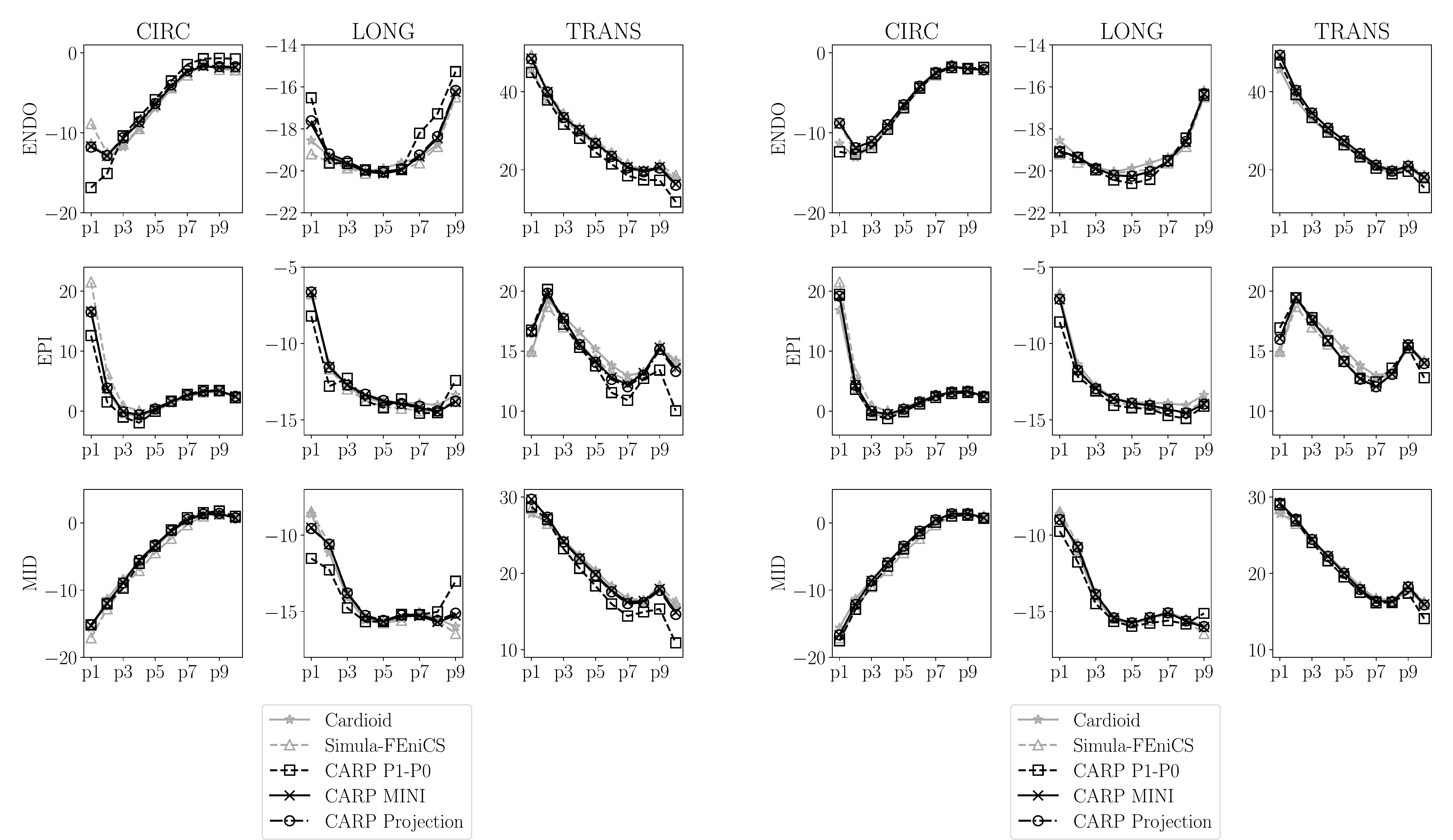}\\[-2em]
  \hspace*{1em}
  {\small {(a)} \emph{$\ell=2$}}
  \hspace*{20.8em}
  {\small {(b)} \emph{$\ell=3$}}\\[-1em]
  \caption{\emph{Ellipsoid benchmark, Guccione material:} longitudinal (LONG),
           circumferential (CIRC), and radial (TRANS) strains at endocardium,
           epicardium, and midwall. Index of points increases from the
         apex to the base. Own results (in black) are compared to benchmark results
         (in gray) presented in~\cite{land2015verification}.}%
  \label{fig:ellipsoid_strain_guccione}
\end{figure*}

\paragraph{Results}
For the transversely isotropic law~\eqref{eq:guccioneStrainEnergy},
we compared our results to selected reference solutions from the benchmark paper~\cite{land2015verification}, namely,
the result from IBM with the Cardioid framework~\cite{Gurev2015high} using P2-P1 elements
and the result from Simula with FEniCS~\cite{logg2012automated} using two-dimensional
P2-P1 elements.
First, the final location of the apex is measured and,
second, circumferential, longitudinal, and radial strains at the endocardium, epicardium,
 and midwall are calculated on points along
apex-to-base lines, see~\cite{land2015verification} for more details.
Results show that the apex location~\Cref{fig:ellipsoid_apex}(a) and strains
\Cref{fig:ellipsoid_strain_guccione} are very similar for the finest level ($\ell=3$) for all
chosen element types. For the level $\ell=2$ the strain solution using simple P1-P0 elements is not
converged showing differences to the benchmark solutions especially in boundary regions at the
apex (p1) and the base (p10), see~\Cref{fig:ellipsoid_strain_guccione}(a).

\begin{figure*}[htbp]
  \includegraphics[width=1.0\linewidth]{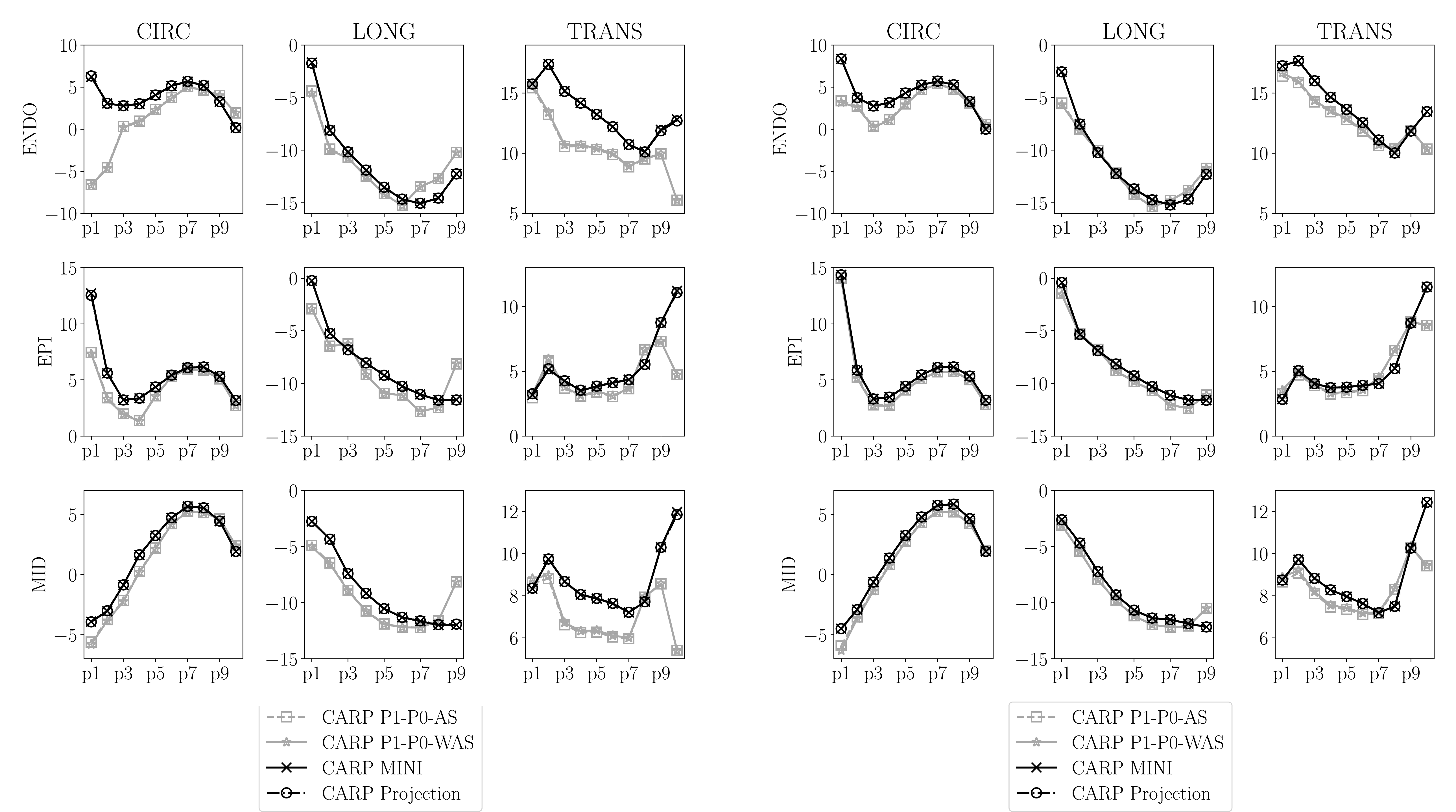}\\[-1.6em]
  \hspace*{1em}
  {\small {(a)} \emph{$\ell=2$}}
  \hspace*{20.8em}
  {\small {(b)} \emph{$\ell=3$}}\\[-1em]
  \caption{\emph{Ellipsoid benchmark, Holzapfel--Ogden material:}
           longitudinal (LONG),
           circumferential (CIRC), and radial (TRANS) strains at endocardium,
           epicardium, and midwall. Index of points increases from the
         apex to the base. P1-P0 elements (in gray) are compared to
         stabilized elements (in black).}%
  \label{fig:ellipsoid_strain_gho}
\end{figure*}

We repeated simulations as above measuring the final apex location and
calculating strains along apex-to-base lines using the orthotropic law~\eqref{eq:benchmark2_ortho}.
We compared results using P1-P0-WAS, P1-P0-AS, projection-stabilized, and MINI elements
in~\Cref{fig:ellipsoid_apex}(b) and~\Cref{fig:ellipsoid_strain_gho}.
Apex locations are very similar for all element types, however, strains
are different, especially in boundary regions close to the apex and the base.
We can see in~\Cref{fig:ellipsoid_strain_gho} that even for the finest
level ($\ell=3$) the strain solution for both P1-P0 formulations is not converged
while solutions for stabilized elements are already very similar for levels
$\ell=2$ and $\ell=3$.


\paragraph{Discussion}
For the transversely-isotropic Guccione material model strains with the presented
projection-stabilized and MINI elements match results using higher order
P2-P1 elements even on coarser grids. Here, also
linear tetrahedral elements seem to be accurate given a fine enough discretization;
this behavior was also observed in~\cite{land2015verification}.

In contrast to that, for the orthotropic Holzapfel--Ogden material, we see in
\Cref{fig:ellipsoid_strain_gho} that P1-P0 elements cannot always accurately
reproduce strains even on the finest grid. On the other hand we can assume that
strains are accurate and almost converged for projection-stabilized and MINI elements as results for levels 2 and 3 are very similar.
Interestingly, for this benchmark, we see no difference in accuracy between the
standard P1-P0 and the P1-P0-WAS formulation.
Most likely the reason for this is
that parameters fitted to human myocardial data~\cite{Guan2019AIC}
are not as stiff in fiber direction compared to the more extreme artificial
benchmark case in~\Cref{sec:holzapfel_test}.
Overall, both approaches with simple linear elements fail to match results
from gold-standard elements, especially in boundary regions.

\subsection{3D-0D closed-loop model of the heart and circulation}%
\label{sec:3D0D}
\paragraph{Simulation setup}
Finally, we show the applicability of our method to an
advanced model of computational cardiac EM.
Here, a 3D model of bi-ventricular EM is coupled to the physiologically
comprehensive 0D CircAdapt model representing atrial mechanics and
closed-loop circulation.
In the present paper, the myocardium of the ventricles was modeled as a nonlinear hyperelastic, (nearly) incompressible and orthotropic material as in~\Cref{eq:HolzapfelAS}.
In particular, for this application, we chose the model proposed
by~\citet{Gultekin2016}
\begin{align*}
  \Theta(J) &:= \ln(J),\qquad
  \overline{\Psi}_{\mathrm{iso}}(\isotensor{C}) := \frac{a}{2b}\left\{
      \exp\left[b(\overline{I}_1-3)\right]-1\right\} \\
  \overline{\Psi}_{\mathrm{aniso}}(\isotensor{C},\vec{f}_0, \vec{s}_0) &:=
    \sum_{i=\mathrm{f,s}}\frac{a_i}{2b_i}\left\{\exp\left[b_i
      (\overline{I}_{4i}-1)^2\right]-1\right\}\nonumber
    + \frac{a_\mathrm{fs}}{2b_\mathrm{fs}}\left\{\exp\left[b_\mathrm{fs}
      (\overline{I}_\mathrm{8fs})^2\right]-1\right\},
\end{align*}
with modified unimodular fourth-invariants to support dispersion of fibers
\[
{I}^\ast_{4i}=\kappa_i\overline{I}_1+(1-3\kappa_i)I_{4i}, \quad i\in\mathrm{f,s}
\]
and standard invariants
\[
  \overline{I}_1 := \mathrm{tr}(\isotensor C),\quad
  \overline{I}_{8\mathrm{fs}}= \vec{f}_0\cdot\isotensor{C}\vec{s}_0.
\]
Analogously, we used the constitutive equation above with
$\Psi_{\mathrm{aniso}}(\tensor{C},\vec{f}_0, \vec{s}_0)$ for the $\mathrm{WAS}$
formulation.

A reaction-eikonal model~\cite{neic2017efficient} was used to generate electrical activation
sequences. Cellular dynamics were described by the Grandi--Pasqualini--Bers model~\cite{Grandi2010a}
coupled to the Land--Niederer model~\cite{Land2012a} to account for length and velocity dependence
of active stress generation, see also~\citet{augustin2016anatomically} for more details on this
strong coupling as well as~\citet{Regazzoni2021} for implementation details on the velocity-dependent
active stress model.
The active stress tensor is computed according to~\cite{Eriksson2013a} as
\[
  \tensor{S}_{\mathrm{a}}
  = S_{\mathrm{a}}\left(\frac{\kappa_\mathrm{f}}{1-2\kappa_\mathrm{f}}\tensor{C}^{-1}
  + \frac{1-3\kappa_\mathrm{f}}{1-2\kappa_\mathrm{f}}(\vec{f}_{0}
  \cdot \tensor{C}\vec{f}_{0})^{-1}\vec{f}_{0}\otimes \vec{f}_{0}\right),
\]
where $S_{\mathrm{a}}$ is the scalar valued active stress generated in the cardiac
myocytes and $\kappa_\mathrm{f}$ is the same dispersion parameter as above.

For the time integration of Cauchy's equation of motion we used a variant of the
generalized-$\alpha$ integrator~\cite{kadapa2017} with spectral radius $\rho_\infty=0$ and
damping parameters $\beta_\mathrm{mass}=\SI{0.1}{\ms^{-1}}$, $\beta_\mathrm{stiff}=\SI{0.1}{\milli \second}$.

\paragraph{Bi-ventricular finite element models}
The bi-ventricular geometry was created according
to~\cite{augustin2021computationally} with an average spatial resolution
of \SI{1.3}{\mm} for the LV and \SI{1.2}{\mm} for the RV.
The resulting mesh used for simulations consisted of \num{557316} elements and \num{111234} nodes.
Fiber and sheet directions were computed by a rule-based method~\cite{bayer2012a:fibers}
with fiber angles changing linearly from $-60^{\circ}$ at the epicardium
to $+60^{\circ}$ at the endocardium~\cite{streeter1969fiber}.

\paragraph{Boundary conditions}
Boundary conditions on the epicardium were modeled using spatially varying
normal Robin boundary conditions~\cite{Strocchi2020pericardium}
to simulate \emph{in-vivo} constrains imposed by the pericardium.
The basal cut plane was constrained by omni-directional spring type boundary conditions.
The 3D ventricular PDE model was coupled to the 0D ODE model \textit{CircAdapt}~\cite{walmsley2015fast}
representing cardiovascular system dynamics according to~\citet{augustin2021computationally}.
An \emph{ex-vivo} setting without pericardial boundary conditions is used to calibrate parameters replicating \emph{ex-vivo} passive inflation experiments
by~\citet{Klotz2007Computational}.
\paragraph{Parameterization}
Passive material parameters
$a=0.4$, $b=6.55$, $a_{\mathrm{f}}=3.05$, $b_{\mathrm{f}}=29.05$,
$a_{\mathrm{s}}=1.25$, $b_{\mathrm{s}}=36.65$,
$a_{\mathrm{fs}}=0.15$, and $b_{\mathrm{fs}}=6.28$
were taken from~\cite{Gultekin2016}.
Dispersion parameters have been identified previously by mechanical experiments
on passive cardiac tissue by~\citet{Sommer2015biomechanical}
and are set to $\kappa_\mathrm{f}=0.08$ and $\kappa_\mathrm{s}=0.09$.

\begin{figure}[hpb]
   \centering
   \begin{subfigure}{0.32\textwidth}
       \includegraphics[width=\textwidth]{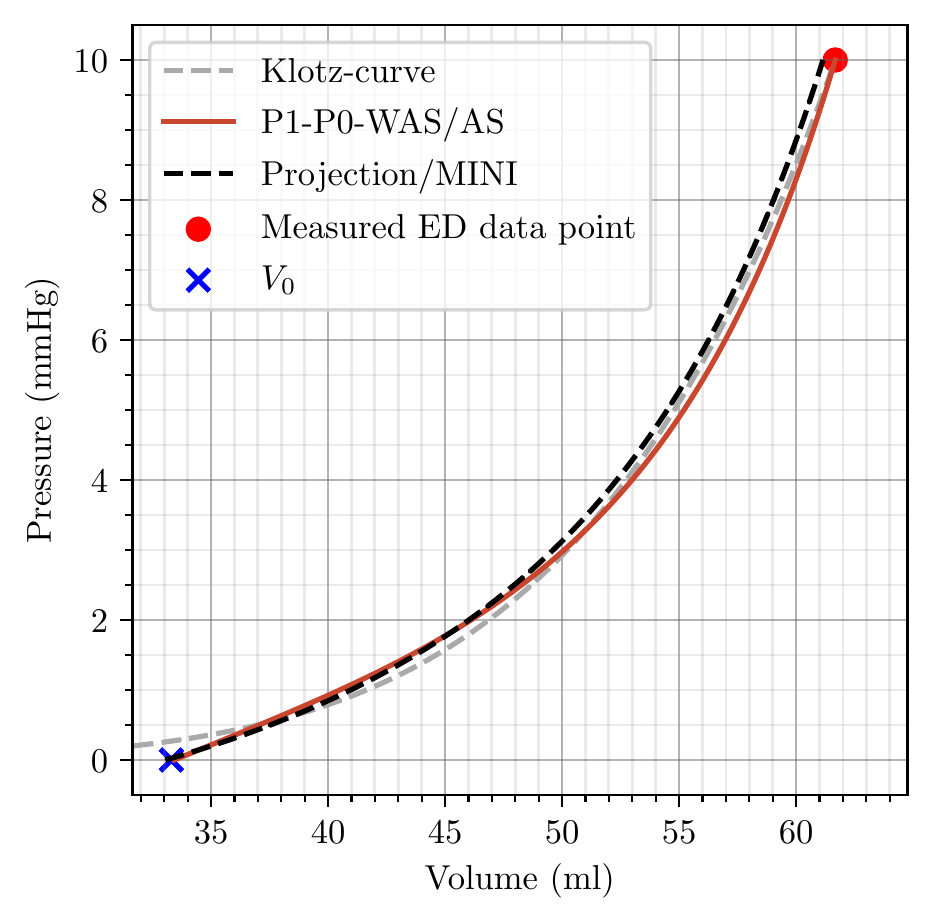}
  \subcaption{}
  \end{subfigure}
   \begin{subfigure}{0.32\textwidth}
       \includegraphics[width=\textwidth]{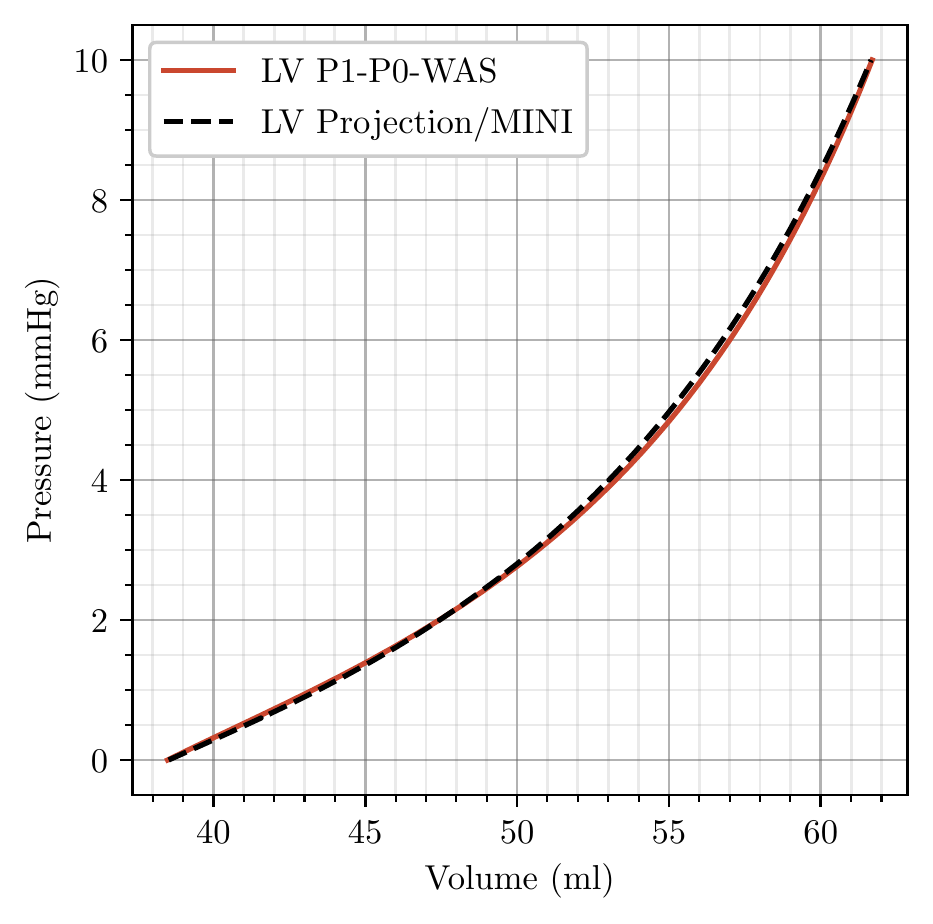}
       \subcaption{}
   \end{subfigure}
   \begin{subfigure}{0.34\textwidth}
     \includegraphics[width=\textwidth]{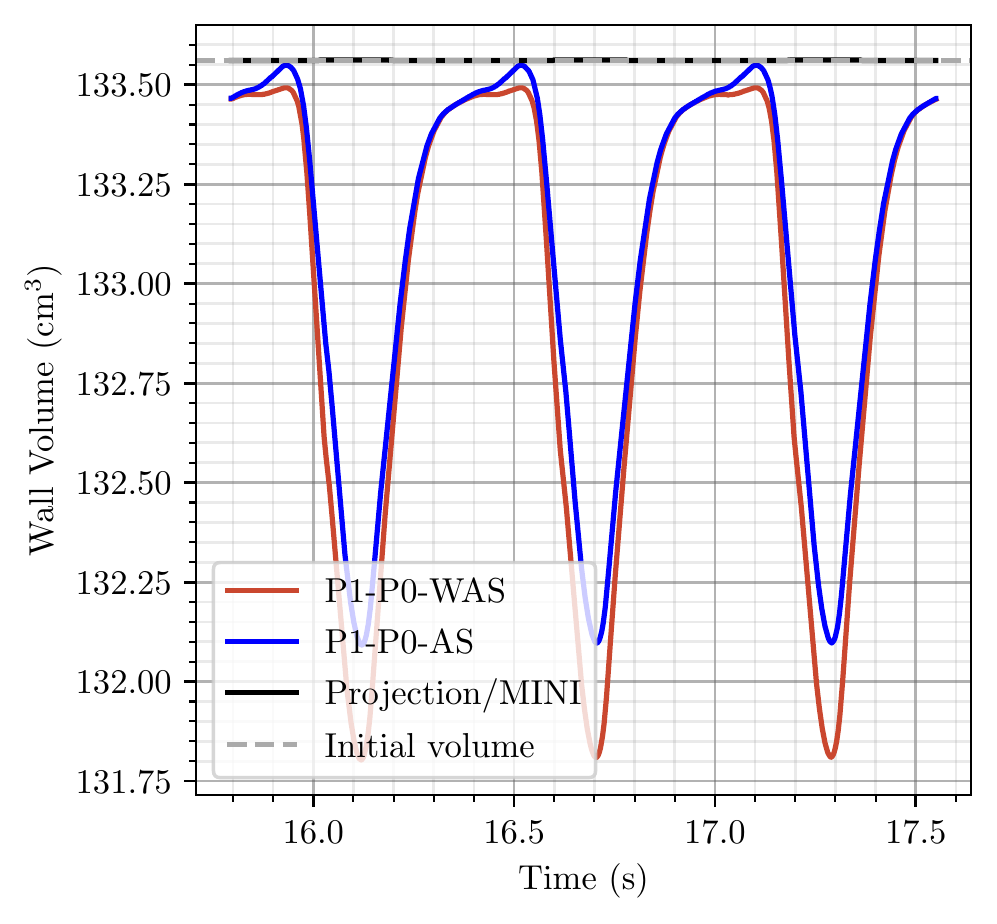}
       \subcaption{}
   \end{subfigure}
       \caption{\emph{3D-0D model of the heart:}
         Passive model calibration and unloading following~\citet{marx2021efficient}.
         (a) An ex-vivo setup was used to calibrate the passive
             material parameters to fit the EDPVR of~\citet{Klotz2007Computational};
         (b) Calibrated material parameters were used in an in-vivo setup
             to generate prestress using a backward-displacement scheme.
         (c) Change in tissue volume over the last three beats.}%
       \label{fig:3D0Dloading}
\end{figure}

To eliminate potential differences in stress/strain
due to parameterization we fitted parameters to achieve similar pressure-volume (PV)
loops and a similar end-diastolic PV relationship (EDPVR) for all element types.
First, initial passive material parameters above were fitted to the empiric
description of the EDPVR by~\citet{Klotz2007Computational}. For each element
type we used a backward displacement algorithm and boundary conditions
replicating experiments in~\cite{Klotz2007Computational} according to~\citet{marx2021efficient},
see~\Cref{fig:3D0Dloading}.
This fitting resulted in multiplicative scaling factors
of \num{0.4529} (P1-P0 elements) and \num{0.9582} (locking-free elements)
for the stress-like material parameters
($a$, $a_{\mathrm{f}}$, $a_{\mathrm{s}}$, $a_{\mathrm{fs}}$);
and in multiplicative scaling factors
of \num{1.0322} (P1-P0 elements) and \num{0.7981} (locking-free elements)
for the dimensionless parameters
($b$, $b_{\mathrm{f}}$, $b_{\mathrm{s}}$, $b_{\mathrm{fs}}$).
The list of fitted passive parameters is given in~\Cref{tab:matparam}.

\begin{figure}[hpb]
   \centering
   \begin{subfigure}{0.42\textwidth}
       \includegraphics[width=\textwidth]{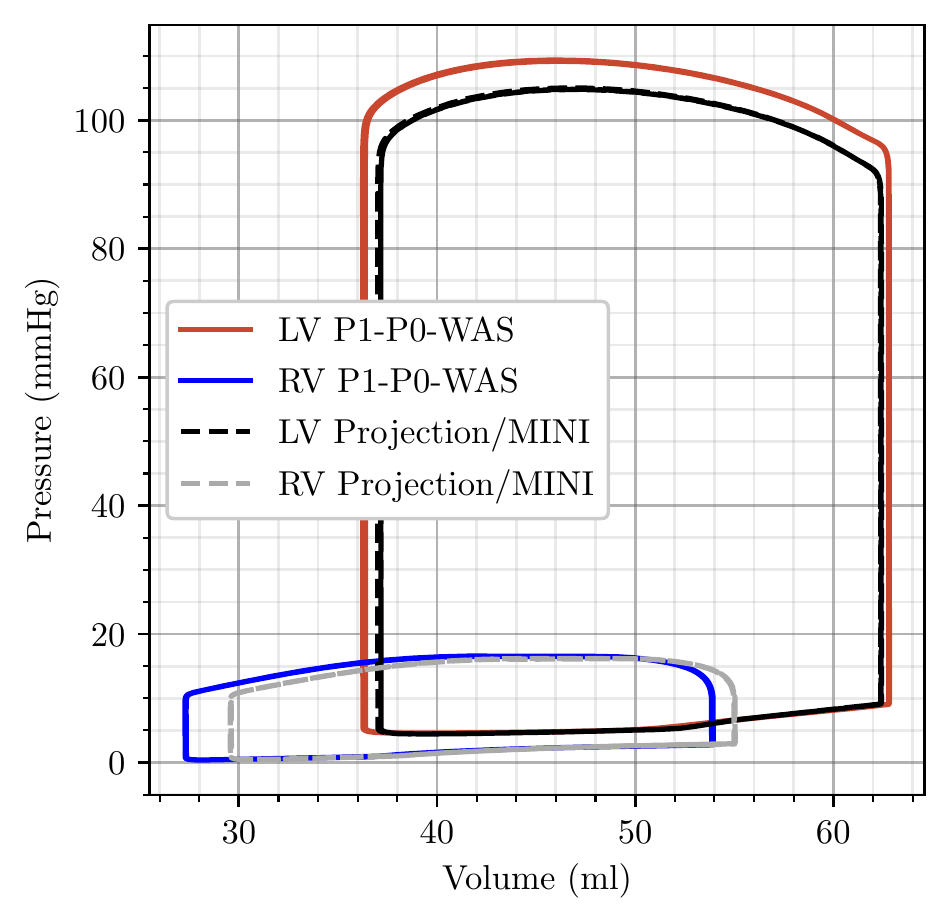}
  \subcaption{}
  \end{subfigure}
   \begin{subfigure}{0.42\textwidth}
       \includegraphics[width=\textwidth]{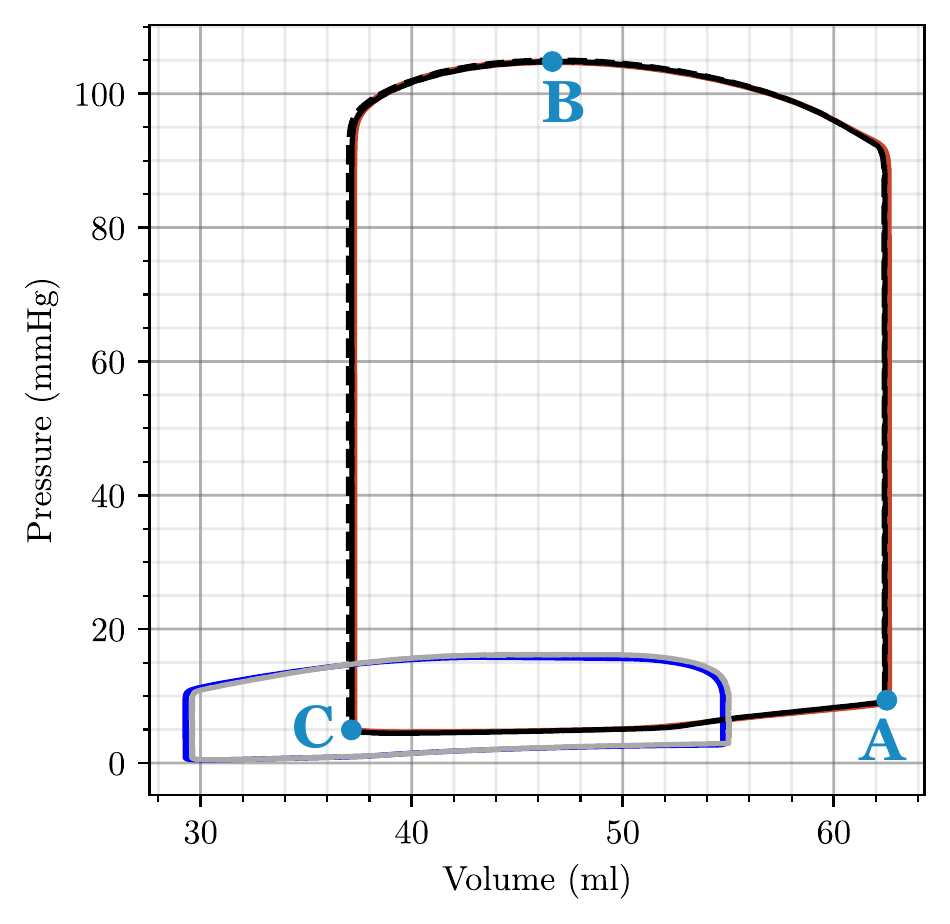}
  \subcaption{}
  \end{subfigure}
   \begin{subfigure}{0.42\textwidth}
       \includegraphics[width=\textwidth]{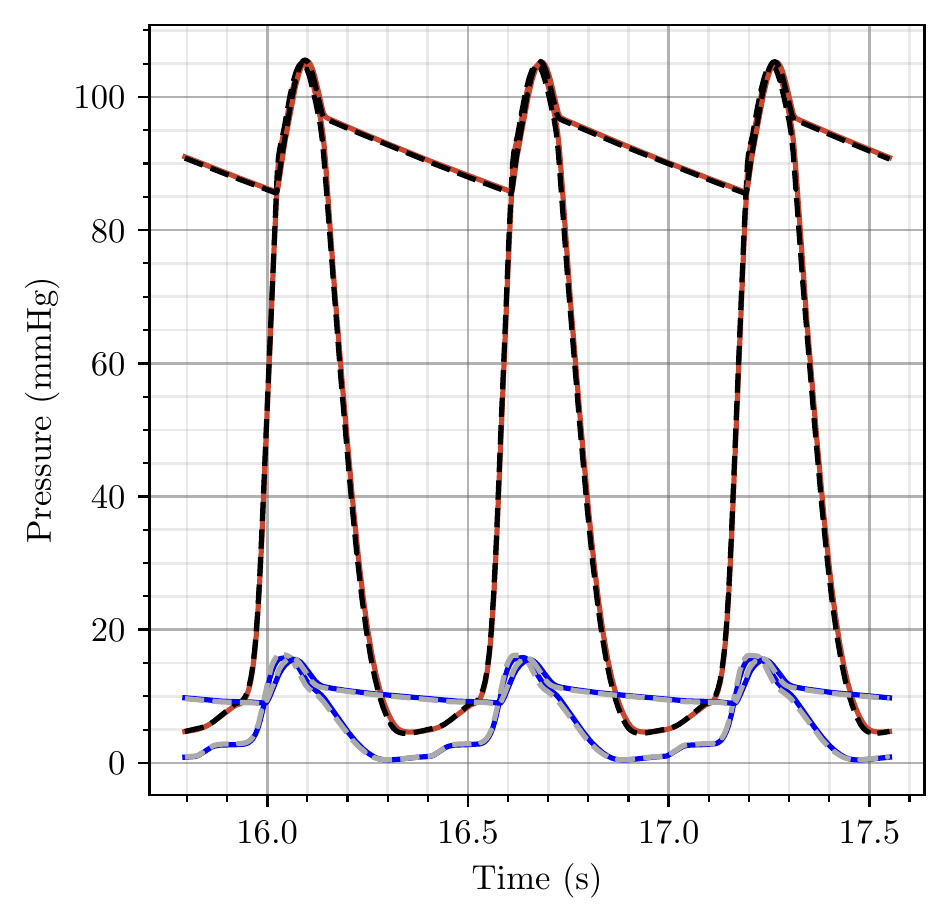}
       \subcaption{}
   \end{subfigure}
   \begin{subfigure}{0.42\textwidth}
       \includegraphics[width=\textwidth]{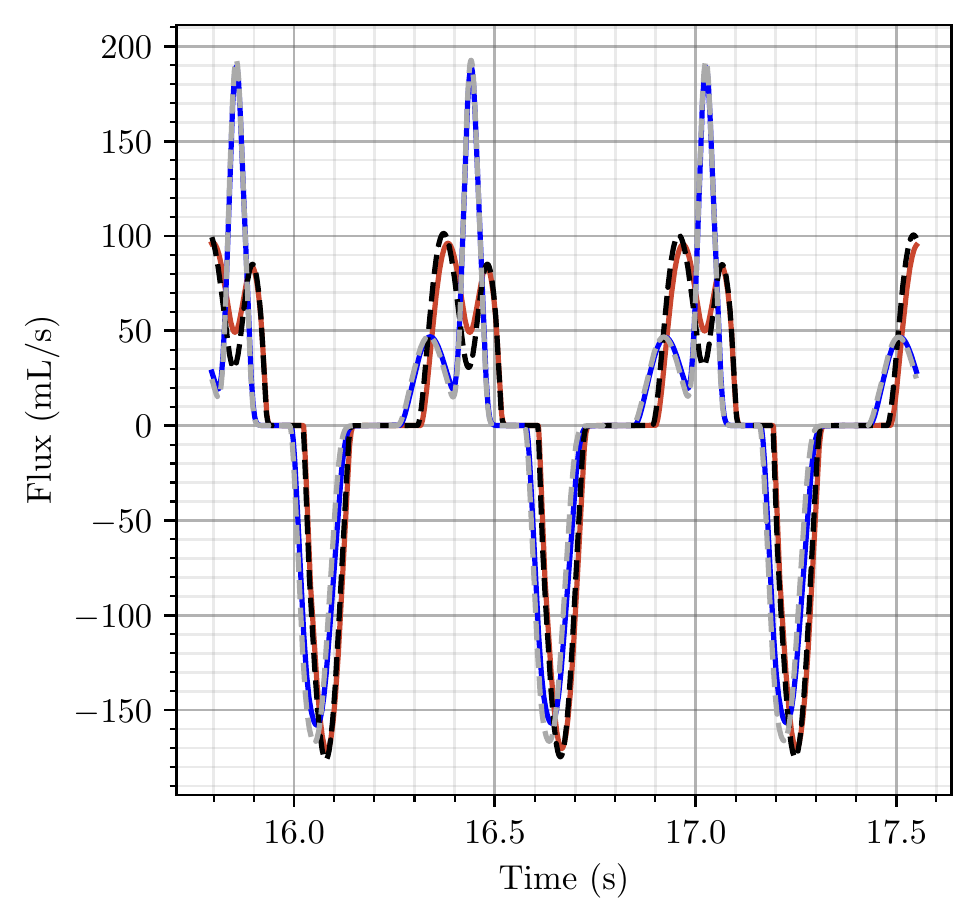}
       \subcaption{}
   \end{subfigure}
       \caption{\emph{3D-0D model of the heart:}
         Plot comparing data traces for the P1-P0-WAS formulation
         (LV: red lines, RV: blue lines) and for the
         MINI and stabilized P1-P1 elements
         (LV: black dashed lines, RV: gray dashed lines).
         Shown are the last 3 beats of the simulation with
         (a) same active stress and preload parameters for P1-P0-WAS
         and locking free elements and
         (b-d) modified active stress and preload parameters for P1-P0-WAS elements
         to reach similar PV loops. In particular we show:
         (a,b) converged PV-loops for both ventricles;
         (c) pressure trace for the LV and RV and pressure in the respective outflow
             vessel;
         (d) in- (negative values) and outflow (positive values) traces of
             both ventricles.
         A, B, C mark the time-points for stress plots
         in~\Cref{fig:3D0Dstresses,fig:3D0Dviolin}.}%
       \label{fig:3D0Dtraces}
\end{figure}

Active stress parameters were fitted using locking-free elements to reach a target peak
pressure of \SI{105}{\mmHg} in the LV.
Using the same active stress parameters, the simulations with P1-P0 elements resulted
in slightly higher peak pressure, \SI{109.2}{\mmHg}, see~\Cref{fig:3D0Dtraces}\;(a).
The end-diastolic state remained almost unchanged while the ejection fraction
increased slightly. We attribute this to a higher contractility of the elements
when the tissue is not modeled as a fully incompressible continuum.
To achieve similar PV loops for P1-P0 elements, simulations were repeated
with reduced active tension $T_\mathrm{ref}^{\bullet}$.
See~\Cref{tab:matparam} for a summary of all active stress parameters.
Simulations parameters of the circulatory system were set as
in~\citet{augustin2021computationally} with a cycle length of \SI{0.585}{\s}.
\paragraph{Results}
First, simulation results in~\Cref{fig:3D0Dloading}\;(a) show that the passive
parameterization - that was performed individually for all element types - allowed
to reach the predicted stress-free volume and the given end-diastolic volume
almost perfectly while reproducing the shape of the Klotz EDPVR curve.
Boundary conditions for this experiment correspond to the~\emph{ex-vivo} setting
described above replicating passive inflation experiments~\cite{Klotz2007Computational}. With fitted material parameters we repeat the backward displacement
algorithm to get a stress-free reference geometry with~\emph{in-vivo}
boundary conditions to model the constrains imposed by the pericardium.
Loading the found stress-free configuration to end-diastolic pressure results
in loading curves as shown in~\Cref{fig:3D0Dloading}\;(b). These loading
curves are almost identical for all element types.
The pre-stressed configuration after this initial loading phase matches the
geometry obtained from imaging and serves as the starting point for
EM heart beat experiments as described in the following.

To get to a converged solution of the closed-loop 3D-0D system we simulated
30 heart beats for each finite element setting: 18 init beats with 1 Newton step
corresponding to a semi-implicit (linearly-implicit) discretization
method~\cite{deuflhard2011newton}; 10 beats with 2 Newton steps which is required
to get a correct update due to the velocity dependence of the active-stress
model; and two final beats with a fully converged Newton
with a relative error reduction of the residual of $\epsilon=10^{-6}$.

See~\Cref{fig:3D0Dtraces}\;(a) for a comparison of the final three pressure-volume loops
with the same active stress parameters which resulted in higher pressures for
P1-P0 elements. In~\Cref{fig:3D0Dtraces}\;(b)--(d) we show traces for the
refitted active stress parameters as described above.
Here, the final three P-V loops in~\Cref{fig:3D0Dtraces}\;(b) coincide.
Hence, the solution is converged and there is also no difference between
the simulation with two Newton steps and the fully converged Newton method.
This also holds for true for pressures in the ventricles and adjacent
arteries, see~\Cref{fig:3D0Dtraces}\;(c), in- and outflow traces~\Cref{fig:3D0Dtraces}\;(d),
and strains/stresses.

Looking at myocardial mass in~\Cref{fig:3D0Dloading}\;(c) we can see that
for projection-stabilized and MINI
elements the mass stays at the initial value of \SI{133.56}{\cm^3} all the time
during the final three beats; this is expected as the tissue is modeled to be
fully incompressible ($1/\kappa=0$).
For P1-P0 elements using a penalty formulation ($\kappa=\SI{650}{\kPa}$)
the tissue is nearly incompressible and especially during the ejection phase the
myocardial mass decreases slightly: maximal \SI{1.17}{\percent} for P1-P0-AS and
\SI{1.38}{\percent} for P1-P0-WAS elements.

\begin{figure}[hpb]
 \centering
   \includegraphics[width=\textwidth]{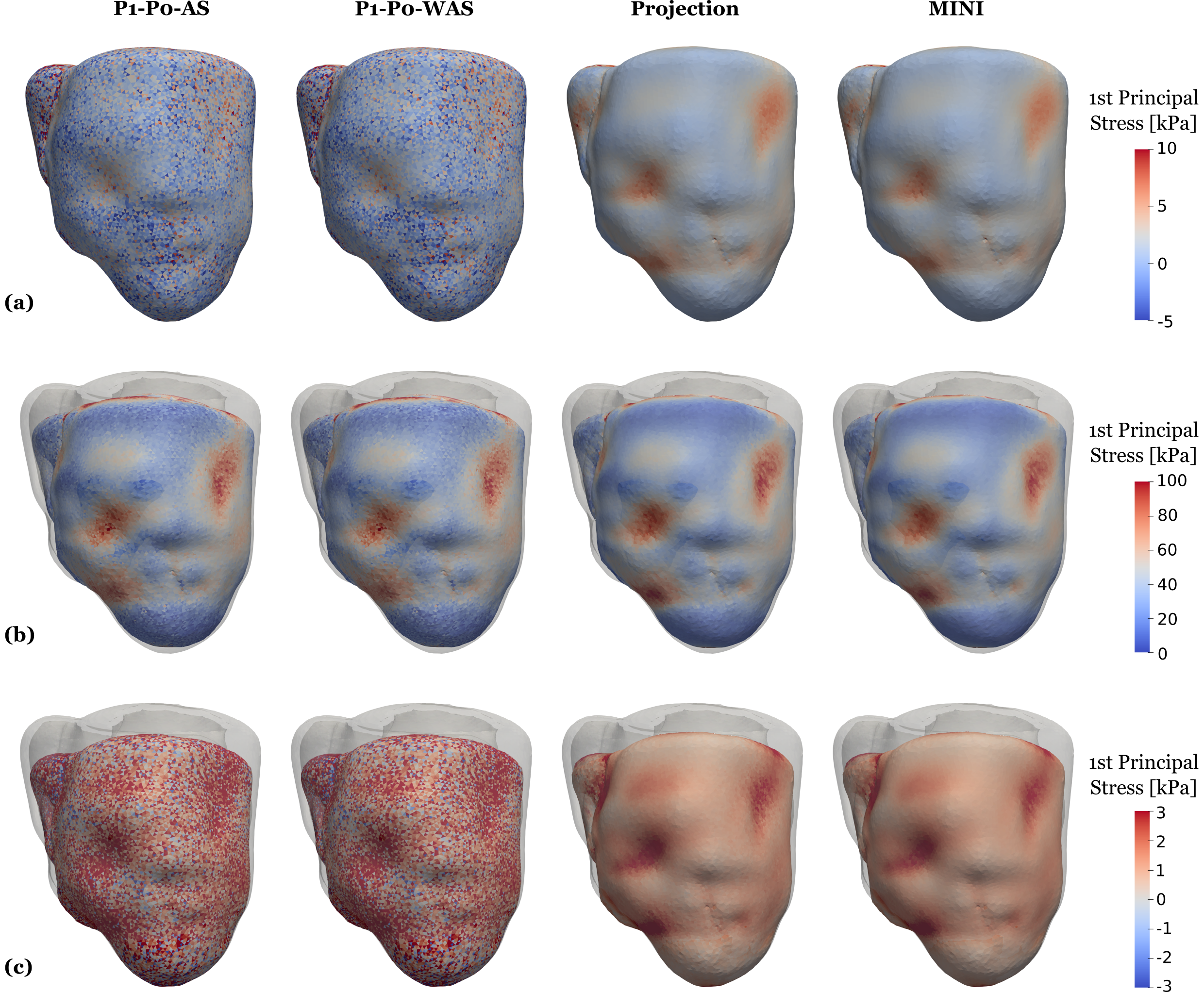}
   \caption{\emph{3D-0D model of the heart:}
     Snapshots of first principal stress values at end-diastole (first row),
     peak systole (second row), and at mitral valve opening (third row),
     see~\Cref{fig:3D0Dtraces} A, B, C for a visualization of the time-points.
     Gray outlines show the end-diastolic configuration.
     Compared are P1-P0-WAS elements (first column),
     stabilized P1-P1 elements (second columns), and MINI elements (third column).
     Shown is total stress as defined tress fields are element-wise and not smoothed.
   }%
   \label{fig:3D0Dstresses}
\end{figure}

\begin{figure}[hpb]
   \centering
   \begin{subfigure}{0.49\textwidth}
       \includegraphics[width=\textwidth]{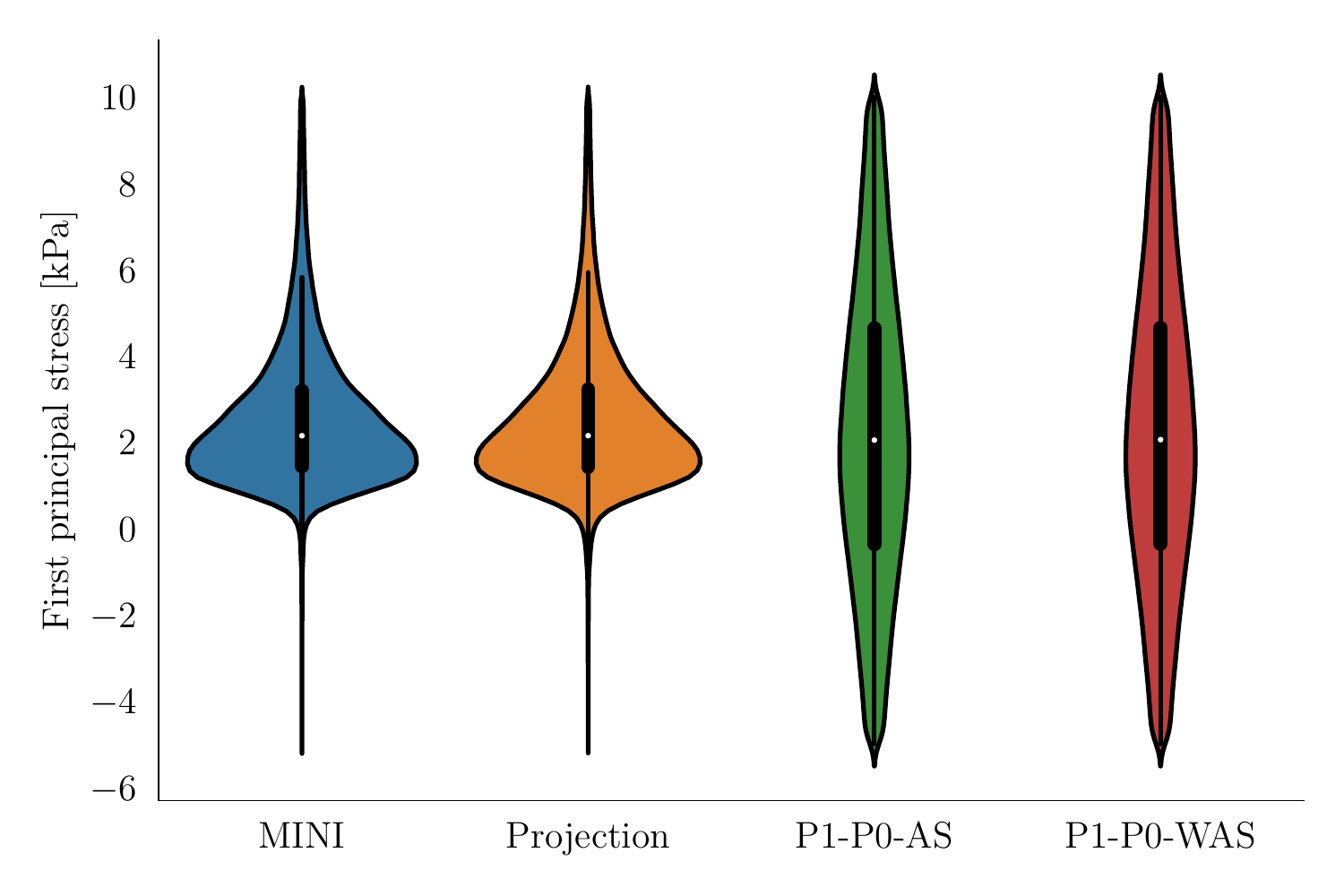}
  \subcaption{}
  \end{subfigure}
   \begin{subfigure}{0.49\textwidth}
       \includegraphics[width=\textwidth]{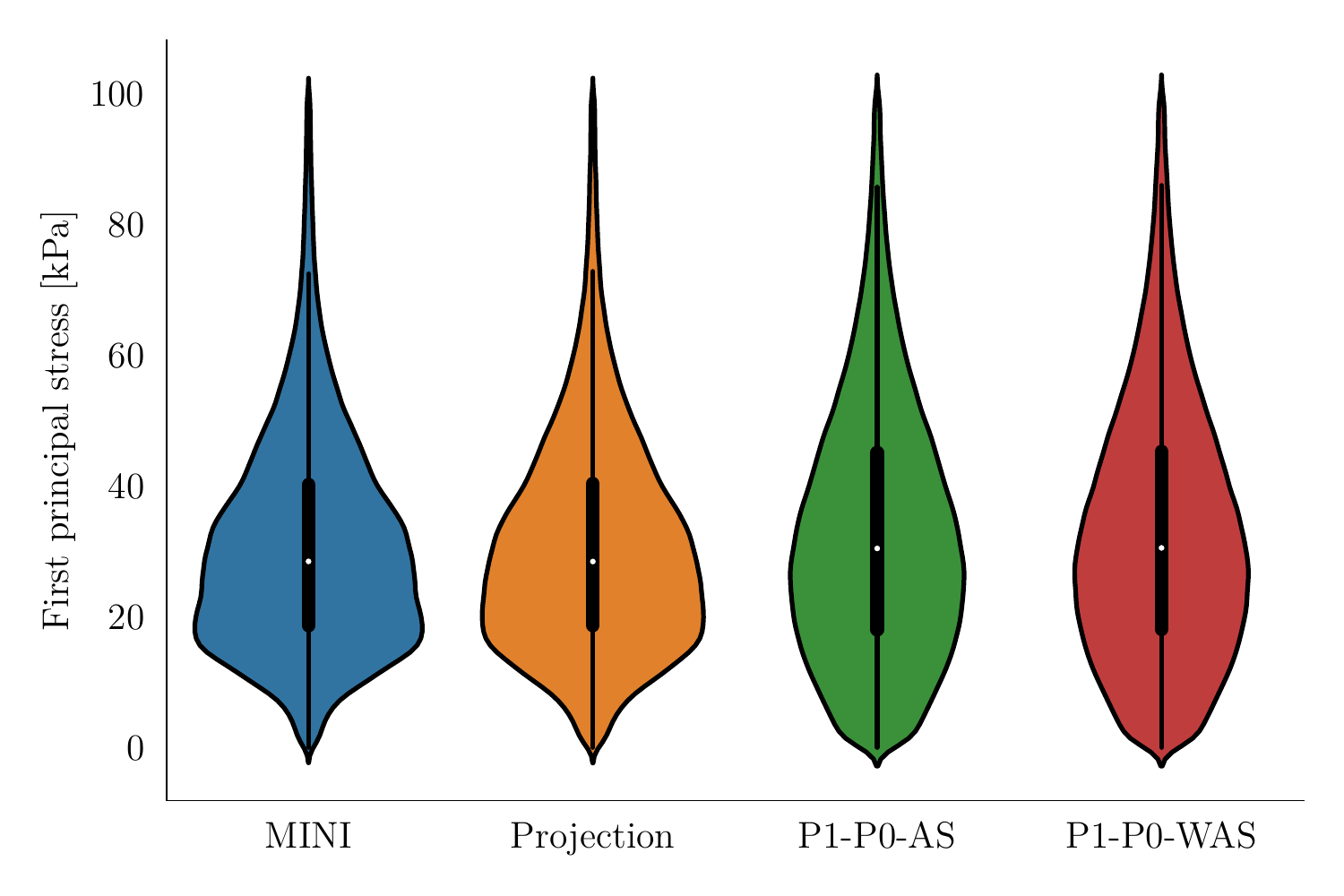}
       \subcaption{}
   \end{subfigure}
   \begin{subfigure}{0.49\textwidth}
       \includegraphics[width=\textwidth]{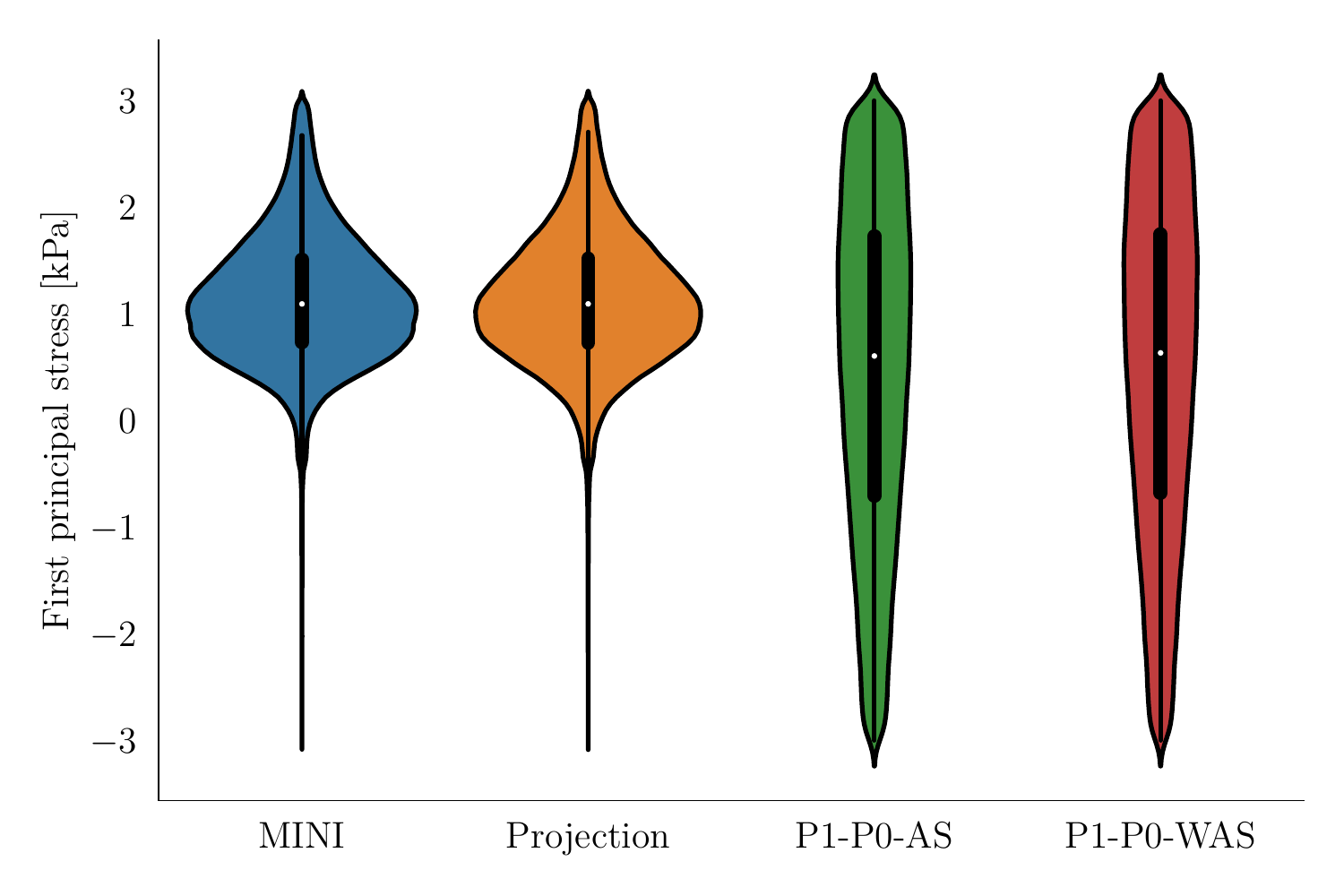}
       \subcaption{}
   \end{subfigure}
   \caption{\emph{3D-0D model of the heart:} Violin plots of the first principal
           stress distribution at (a) end-diastole, (b) peak systole,
           (c) mitral valve opening.}%
       \label{fig:3D0Dviolin}
\end{figure}

In contrast to pressure, volume, and flow traces, stresses show a very different
pattern when comparing locking-free to simple P1-P0 elements,
see~\Cref{fig:3D0Dstresses}. In this plot we show element-wise, total first
principal stress at three time points marked by A, B, and C in the P-V loop
in~\Cref{fig:3D0Dtraces}\;(b) which represent A, the most expanded (end-diastole),
B, the highest total stress (peak systole),
and C, the most contracted (beginning of filling phase)
states of the ventricles.
Especially at end-diastole~\Cref{fig:3D0Dstresses}\;(a) and the beginning of
the filling phase~\Cref{fig:3D0Dstresses}\;(c) where passive stress dominates
and active stress is close to zero we see a distinct checkerboard pattern
for P1-P0 elements while solutions for projection-stabilized and MINI elements
are smooth. Also in violin plots showing the stress distribution over the
whole tissue domain we see a clear difference for these time points,
see~\Cref{fig:3D0Dviolin}\;(a) and (c).
On the other hand, the stress distribution is very similar for all element
types at peak-systole where active stress dominates,
see~\Cref{fig:3D0Dstresses}\;(b) and~\Cref{fig:3D0Dviolin}\;(b).

\begin{table}
\caption{Summary of electrical and mechanical material parameters.}%
\small\centering
\label{tab:matparam}
  \begin{tabular}{r@{\hskip 0.8ex}lr@{\hskip 0.8ex}lr@{\hskip 0.8ex}lr@{\hskip 0.8ex}l}
\toprule
  \multicolumn{8}{l}{\emph{Passive stress parameters: P1-P0 elements}} \\
  $a=$ & \SI{0.1812}{\kPa}, & $a_{\mathrm{f}}=$&\SI{1.3813}{\kPa}, &
  $a_{\mathrm{s}}=$& \SI{0.5661}{\kPa}, & $a_{\mathrm{fs}}=$&\SI{0.0679}{\kPa}, \\
  $b=$ &\SI{6.7609}{[-]},  & $b_{\mathrm{f}}=$ & \SI{29.9854}{[-]}, &
  $b_{\mathrm{s}}=$&\SI{37.8301}{[-]}, & $b_{\mathrm{fs}}=$&\SI{6.4822}{[-]}. \\
  \midrule
  \multicolumn{8}{l}{\emph{Passive stress parameters: locking-free elements}} \\
  $a=$ & \SI{0.3833}{\kPa}, & $a_{\mathrm{f}}=$&\SI{2.9225}{\kPa}, &
  $a_{\mathrm{s}}=$& \SI{1.1978}{\kPa}, & $a_{\mathrm{fs}}=$&\SI{0.1437}{\kPa}, \\
  $b=$ &\SI{5.2278}{[-]},  & $b_{\mathrm{f}}=$ & \SI{23.1848}{[-]}, &
  $b_{\mathrm{s}}=$&\SI{29.2504}{[-]}, & $b_{\mathrm{fs}}=$&\SI{5.0121}{[-]}. \\
  \midrule
  \multicolumn{8}{l}{\emph{Active stress parameters}} \\
    $T^\mathrm{LV}_{\mathrm{ref}}=$&$\SI{200.0}{\milli\N\per\square\mm}$,  &
    $T^\mathrm{RV}_{\mathrm{ref}}=$&$\SI{160.0}{\milli\N\per\square\mm}$,  &
  ${[\textrm{Ca}^{2+}]}_{\mathrm{T50}}=$&\SI{0.52}{\micro\mole\per\liter}, &
    $\mathrm{TRPN}_{50}=$&\SI{0.37}{[-]}, \\
    $n_{\mathrm{TRPN}}=$&\SI{1.54}{[-]}, & $k_{\mathrm{TRPN}}=$&$\SI{0.14}{\per\milli\s}$,  &
    $n_{\mathrm{xb}}=$&\SI{3.38}{[-]},  & $k_{\mathrm{xb}} =$&$ \SI{4.9e-3}{\per\milli\s}$. \\
  \midrule
  \multicolumn{8}{l}{\emph{Adapted active stress parameters for P1-P0 elements}} \\
    $T^\mathrm{LV}_{\mathrm{ref}}=$&$\SI{190.0}{\milli\N\per\square\mm}$,  &
    $T^\mathrm{RV}_{\mathrm{ref}}=$&$\SI{130.0}{\milli\N\per\square\mm}$.  &&\\
  \bottomrule
  \end{tabular}
\end{table}

\paragraph{Numerical performance} Computational times for the simulation using
different element types are given in~\Cref{tab:num_perf}.
Simulations were performed on \num{128} cores of Archer2
and we distinguish between solver-time, $t_{\rm{s},\bullet}$,
the accumulated time of the linear solver (GMRES),
and assembly-time, $t_{\rm{a},\bullet}$, the accumulated time of matrix and vector assembling of the linearized system~\eqref{eq:linsys1}--\eqref{eq:linsys3}.

In total, for a full simulation with loading, 18 initialization beats with 1
Newton step, 10 initialization beats with 2 Newton steps,
and 2 final beats with a fully converging Newton method the computational
costs were around \SI{3.5}{\hour} for P1-P0 elements,
\SI{13}{\hour} for projection stabilized elements,
and \SI{17.5}{\hour} for MINI elements,
see~\Cref{tab:num_perf} for exact values.
Here, in addition to GMRES solver and assembly times,
also input-output times, the solution of the R-E model governing EP,
ODE times, and postprocessing are taken into account.
Using a coarser mesh with \num{45686} elements and \num{11850} nodes tractable
computational times could also be achieved on a desktop machine (AMD Ryzen Threadripper 2990X)
with 2.5 minutes for one heart beat on 32 cores using P1-P0 elements and \SI{7.15}{minutes} using
locking-free projection stabilized elements.
Total computational times on the desktop machine for 30 beats were
\SI{94.6}{minutes} for P1-P0 and \SI{264.4}{minutes} for projection stabilized elements.

\begin{figure*}[htbp]
  \centering
  \includegraphics[width=0.85\linewidth,keepaspectratio]{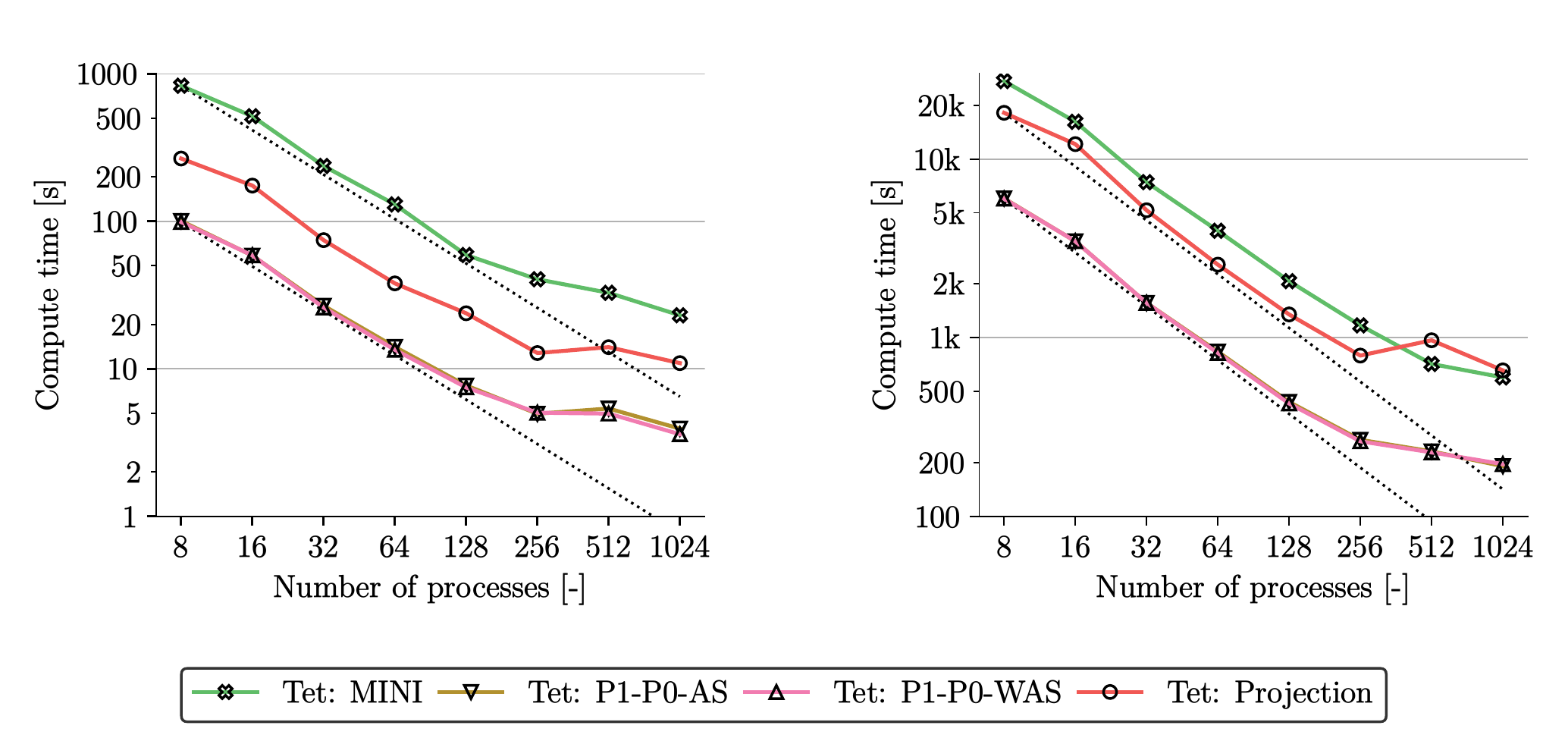}
  \caption{\emph{3D-0D model of the heart:} Strong scaling results for the loading phase (left) and one beat (right).
  Simulations were performed on 16 to 1024 cores on Archer2.}%
  \label{fig:pc:strong_scaling}
\end{figure*}

In~\Cref{fig:pc:strong_scaling} we show strong scaling properties of the
simulation on 16 to 1024 cores of Archer2. Loading and heart beat experiments
scale well up to 256 cores for all element types. For 512 and 1024 cores strong
scaling efficiency drops markedly due to small local partition sizes
($<500$ degrees of freedom per partition).

\begin{table}[htbp]\small\centering
  \begin{tabular}{lrrrrrrr}
  \toprule
    Type & DOF &
    $t_{\rm{s},1}/t_{\rm{a},1}$ &
    $t_{\rm{s},2}/t_{\rm{a},2}$ &
    $t_{\rm{s,c}}/t_{\rm{a,c}}$ &
    $T_{\rm{b,1}}/T_{\rm{b,2}}/T_{\rm{b,c}}$  &
    $T_\mathrm{ld}$ & Total\\
    $[-]$  & $[-]$ & $[\si{\second}]$ & $[\si{\second}]$ &
    $[\si{\second}]$ & $[\si{\second}]$ & $[\si{\second}]$ &
    $[\si{\min}]$ \\
  \midrule
  P1-P0-AS  & \num{333702}
  & \num{ 93.7}/\num{93.4}
  & \num{183.6}/\num{186.8}
  & \num{938.4}/\num{739.0}
  & \num{221.3}/\num{438.6}/\num{1947.5}
  & \num{7.71}& \num{207.2}\\
  P1-P0-WAS & \num{333702}
  & \num{ 92.1}/\num{92.9}
  & \num{181.2}/\num{185.8}
  & \num{922.4}/\num{739.7}
  & \num{218.7}/\num{434.4}/\num{1930.4}
  & \num{7.48} & \num{205.0} \\
 Projection & \num{444936}
            & \num{573.6}/\num{119.1}
            & \num{1289.1}/\num{238.2}
            & \num{7479.6}/\num{921.3}
            & \num{727.1}/\num{1606.0}/\num{8694.3}
            & \num{23.8} & \num{777.0}\\
  MINI      & \num{444936}
            & \num{600.8}/\num{495.2}
            & \num{1242.1}/\num{990.4}
            & \num{6126.3}/\num{3078.8}
            & \num{1181.8}/\num{2328.3}/\num{9547.5}
            & \num{59.1} & \num{1060.5} \\
  \bottomrule
  \end{tabular}
  \caption{Summary of computational times on 128 cores of Archer2 for the
      different finite element types.
    Given are solver, assembly, and total computational times for one beat using
    one Newton iteration ($t_{\rm{s},1}$, $t_{\rm{a},1}$)
    two Newton iterations ($t_{\rm{s},2}$, $t_{\rm{a},2}$)
    and a fully converged Newton solution ($t_{\rm{s},c}$, $t_{\rm{a},c}$).
    $T_{\rm{b,1}}$, $T_{\rm{b,2}}$, and $T_{\rm{b,c}}$ correspond to
    the total simulation time per heart beat for one Newton iteration,
    two Newton iteration, and fully converged Newton scenarios.
    Timings refer to a single heart beat lasting \SI{0.585}{\second}
    at a time step size of \SI{1}{\milli \second}.
    In addition, the times required for the loading phase, $T_\mathrm{ld}$,
    using \num{32} load steps and the total simulation times including
    0D solution, IO, and postprocessing are presented.}%
    \label{tab:num_perf}
\end{table}
\paragraph{Discussion}
In this benchmark, we show the most complete model of cardiac EM that is
currently available:
i) Cardiac electrophysiology was modeled by a reaction-Eikonal model which
predicts potential fields with high fidelity even on coarser
grids~\cite{neic2017efficient}.
ii) Cellular dynamics were modeled by the physiological
Grandi--Pasqualini--Bers model~\cite{Grandi2010a} which is coupled to
the Land--Niederer model~\cite{Land2012a}.
This allows for strong coupling, i.e., to account for length and velocity
effects on the cytosolic calcium transient, using an approach as
described in~\citet{augustin2016anatomically}.
iii) Passive tissue mechanics was modeled by the recent Holzapfel--Ogden type model~\cite{Gultekin2016} and active stress according to~\cite{Eriksson2013a}
using a recent approach by \citet{Regazzoni2021} to avoid oscillations.
Note that both, passive and active stress computation, account for fiber
dispersion, hence, allowing to model the active tension generated by dispersed fibers.
iv) Spatially varying Robin boundary conditions were included to model the effect of the pericardium~\cite{Strocchi2020pericardium}.
v) The 3D PDE model was coupled to the physiologically comprehensive 0D closed-loop
model CircAdapt of the cardiovascular system. This allows to replicate physiological behaviors
under experimental standard protocols altering loading conditions and contractility~\cite{augustin2021computationally}.
vi) Simulations were performed using locking-free finite elements - as presented
in this paper - to accurately compute stress distributions.

All computational models of cardiac EM presented in the literature so far,
e.g., \cite{Gurev2015high,augustin2021computationally,piersanti20213d0d,Kariya2020,Hirschvogel2017monolithic,Pfaller2019,Wang2021,gerach2021electro}, are missing one or mostly more of the above points.
While the importance of model components i) -- v) was discussed extensively
in references above we could show in this paper that also vi) is necessary
to compute accurate stress fields in the tissue. Depending on the
application this could be a very critical modeling component,
e.g., for the accurate prediction of rupture risks or for the estimation of
growth and remodeling based on stress.
Note that also for this benchmark we see no substantial difference between
simulation outcomes using standard standard P1-P0 and the P1-P0-WAS formulation; both approaches fail to match stress results from gold-standard elements.

While stress fields differ vastly the PV loops computed with locking-free and
simple P1-P0 elements are very similar; at least if active and passive tissue
parameters are fitted independently for each element type.
In particular, passive parameters fitted to the Klotz curve using P1-P0
elements correspond to a softer material. Here, the fitting compensates
locking effects to a certain degree.
On the other hand, the reference peak tension parameter ($T_\mathrm{ref}^{\bullet}$)
had to be slightly reduced to reach the same target value as
projection stabilized and MINI elements.
In this case, the softer passive material and the higher compressibility
of the tissue using the penalty term in the P1-P0 formulation are most likely
the reasons for higher active stress generation.
With these adaptations, PV traces for the different element types are
almost identical, see~\Cref{fig:3D0Dtraces}.
This shows that simple P1-P0 elements can predict most
simulation outputs as good as gold standard formulations and are thus
adequate for cardiac EM simulations when stresses are not an outcome or
critical simulation value. In this case, the computational efficiency of
P1-P0 elements might trump the numerical accuracy of locking-free elements.

High resolution EM models require efficient numerical solvers to limit the
computational cost that results from a high number of degrees of freedom to
capture anatomical details as well as a high number of time steps.
Strong scaling characteristics of our EM framework was reported in detail
previously~\cite{augustin2016anatomically,karabelas2018towards}.
In the present work, we showed that strong scaling is preserved when using
locking free elements and a coupling to a 0D model of the circulatory system.

Using the advanced approach presented in this paper, the time needed for
the passive filling of the bi-ventricular model of the heart is very low.
For the projection stabilized element loading times are less than half a
minute on 128 cores@2.25Ghz of Archer2~(\Cref{tab:num_perf}) for the simulation
with \num{444936} degrees of freedom.
Even on 8 cores the passive filling could be achieved within 5
minutes~(\Cref{fig:pc:strong_scaling}).
In comparison, a recent work~\cite{hurtado2021accelerating} reports
compute times for a similar passive inflation scenario using locking-free
elements that were around 162 minutes (\num{203214} degrees of freedom, 8 cores@2.8Ghz).
Fast loading times are crucial for the estimation of the stress-free reference configuration using fixed-point iterations~\cite{Sellier2011,Rausch2017Augmented}.

Computational cost for one heart beat -- using grids with a comparable number
of elements and nodes but, in general, P1-P0 elements --
range from 1.8 to 24 hours in previous
studies~\cite{Gurev2015high,Kariya2020,Hirschvogel2017monolithic,Pfaller2019,Wang2021}.
In this work, we could show that even with locking-free elements one heart beat can be simulated within 27 minutes on 128 cores of Archer2,
see~\Cref{tab:num_perf}, and within 11 minutes using 1024 cores,
see~\Cref{fig:pc:strong_scaling}.
Still, P1-P0 elements are computationally less expensive with one heart beat in around 7 minutes using 128 cores and 3.2 minutes using 1024 cores.
Using a coarser mesh as in~\cite{augustin2021computationally} fast
computational times are also possible on desktop machines with 2.5 minutes
for one heart beat on 32 cores using P1-P0 elements and 7.15 minutes using
locking-free elements.
This computational efficiency is of paramount importance for future
parameterization studies where numerous forward simulations have to be carried
out to personalize models to patient data.

\paragraph{Limitations}
First, we set an arbitrary number of 30 heart beats which was more than enough to
reach a limit cycle in all experiments. However, an automatic stopping criterion
could be used that stops the simulation after reaching the limit cycle.
Simulation times could be further reduced by accelerated the convergence to a
limit cycle using data-driven 0D emulators~\cite{regazzoni2021accelerating}
or by tuning the 0D CircAdapt model to predict P-V traces from the 3D-0D model
in a better fashion, hence, reducing the number of beats needed to a converged
solution.

Second, the parameterization of the model has room for improvement. Passive
parameters were fitted to the empiric Klotz curve using end-diastolic volume
and pressure and active stress was fitted using a target peak pressure value.
Other model components such as the reaction-Eikonal model and CircAdapt were not
parameterized and we were using default values from the literature.
The personalization of the complete model to patient-specific data is not
within the scope of this contribution, however, the computational efficiency
of the model is of crucial importance for parameter identification studies that
often require a large number of forward simulations.

Third, no independent validation of the model is performed. This could be done
by comparing displacements or strains predicted by the model
to observations from cine MRI or 3D tagged MRI data, see, e.g.,
\cite{Ponnaluri2019, Pfaller2019,sack2018construction}.
However, in this work, we focused on showing advantages of locking-free
elements for applied simulations using an advanced setup which is necessary
to replicate physiological behavior. A rigorous, independent validation against
for several patient-specific cases using image data will be the
focus of future studies.

Finally, locking-free formulations as presented in this paper require the
solution of a block system, which in turn necessitates suitable preconditioning for
computational efficiency. This is not a trivial task, however, preconditioners
used for simulations in this paper are publicly available through the open-source
software framework PETSc.

\section{Conclusion}\label{sec:conclusion}
In this study, we introduced stabilization techniques that
accelerate simulations of anisotropic materials, in particular,
nearly and fully incompressible fiber-reinforced solids such as arterial wall
or myocardial tissue.
A MINI element formulation and a simple and computationally efficient technique
based on a local pressure projection were presented.
Both methods were applied for the first time for simulations of anisotropic
materials and showed to be an excellent choice when the use of higher order or
Taylor--Hood elements is not desired. This is the case, e.g.,
for detailed, high-resolution problem domains that results in a high number
of degrees of freedom.
We showed that both approaches are very versatile and can be applied to stationary and
transient problems as well as hexahedral and tetrahedral grids without modifications.
It is worth noting that all required implementations are purely on the element level, thus,
facilitating an inclusion in existing finite element codes.
Furthermore, solvers and preconditioners used to solve the linearized block system
of equations are available through the open-source software package
\emph{PETSc} \cite{petsc-user-ref}.

We showed the robustness and accuracy of the chosen approaches in two benchmark problems
from the literature: first, a thick-walled cylindrical tube representing arterial tissue
and second, an ellipsoid representing LV myocardial tissue.
Additionally, in a third application of the stabilization approaches, we presented a
complex 3D-0D model of the ventricles.
This constitutes the first computational EM model of the heart where
all components are captured by physiological, state-of-the-art models.
We could show that for the first time accurate and physiological
cardiovascular simulations are feasible within a clinically tractable time frame.

Computational efficiency of the methods is unprecedented in the literature and
the framework shows excellent strong scaling on desktop and HPC architectures.
The high versatility of the one-fits-all approach allows the simulation of
nearly and fully incompressible fiber-reinforced materials in many different scenarios.
Overall, this offers the possibility to
perform accurate simulations of biological tissues in clinically tractable
time frames, also enabling parameterization studies where numerous
forward simulations have to be carried out to personalize models to patient data.

%
\section*{Declaration of competing interest}
The authors declare that they have no known competing financial interests or
personal relationships that could have appeared to influence the work reported in this paper.

\section*{Acknowledgements}
This project has received funding from the European Union's Horizon 2020 research and
innovation programme under the Marie Sk{\l}odowska–Curie Action H2020-MSCA-IF-2016 InsiliCardio, GA No. 750835
and under the ERA-NET co-fund action No. 680969 (ERA-CVD SICVALVES, JTC2019) funded by the Austrian Science Fund (FWF), Grant I 4652-B to CMA.
Additionally, the research was supported by the Grants F3210-N18 and I2760-B30 from the Austrian Science Fund (FWF)
and a BioTechMed Graz flagship award ``ILearnHeart'' to GP.
Further, the project has received funding from the European Union's Horizon 2020 research and innovation programme under the ERA-LEARN co-fund
action No. 811171 (PUSHCART, JTC1\_27) funded by ERA-NET ERACoSysMed to GP.

\clearpage
\bibliographystyle{elsarticle-num-names}
\bibliography{plain.bib}

\newgeometry{margin=1in}

\appendix
\section{Linearization}
\label{sec:appendix_linearization}
We will give a short summary of the linearization of a cavity volume $V_\mathrm{CAV}$ defined by
\begin{align*}
  V_{\mathrm{CAV}} := \frac{1}{3} \int\limits_{\Gamma_{\mathrm{CAV}}} \vec{x} \cdot \vec{n} \dsx.
\end{align*}
Using Nanson's formula and $\vec{x} = \Xvec{X} + \vec{u}$ we can rewrite this as
\begin{align*}
 V_{\mathrm{CAV}} = \frac{1}{3} \int\limits_{\Gamma_{\mathrm{CAV},0}} (\Xvec{X} + \vec{u}) \cdot J \tensor{F}^{-\top} \Xvec{N}\dsX
\end{align*}
Using the known linearizations
\begin{align}
  \frac{\partial J}{\partial \tensor{F}} : \Grad \Delta \vec{u} &= J \tensor{F}^{-\top} : \Grad \Delta \vec{u} \\
  \frac{\partial \tensor{F}^{-\top}}{\partial \tensor{F}} : \Grad \Delta \vec{u} &= -\tensor{F}^{-\top}(\Grad \Delta \vec{u})^\top \tensor{F}^{-\top}
\end{align}
we can calculate the linearization around $\Delta \vec{u}$ as
\begin{align}
  \label{eq:var_volume_cav}
 d_k(\Delta \vec{u}) := D_{\Delta \vec{u}} V_{\mathrm{CAV}} &= D_{\Delta \vec{u}} \frac{1}{3}\int\limits_{\Gamma_{\mathrm{CAV}}} \vec{x} \cdot \vec{n} \dsx \\
                                    &= D_{\Delta \vec{u}} \frac{1}{3}\int\limits_{\Gamma_{\mathrm{CAV},0}} J \left(\vec{X} + \vec{u}\right) \cdot \tensor{F}^{-\top} \Xvec{N} \dsX \\
                                    &= \frac{1}{3} \int\limits_{\Gamma_{\mathrm{CAV},0}} J (\tensor{F}^{-\top} : \Grad \Delta \vec{u}) \vec{x} \cdot \tensor{F}^{-\top} \Xvec{N} \dsX\\
                                    & - \frac{1}{3} \int\limits_{\Gamma_{\mathrm{CAV},0}} J \vec{x} \cdot \tensor{F}^{-\top}(\Grad \Delta \vec{u})^\top \tensor{F}^{-\top} \Xvec{N}\dsX \\
                                    &+ \frac{1}{3} \int\limits_{\Gamma_{\mathrm{CAV},0}} J \Delta \vec{u} \cdot \tensor{F}^{-\top}\Xvec{N}\dsX
\end{align}

\section{Static Condensation for Inhomogeneous Neumann Boundary Condition}
\label{sec:appendix_static_condensation}
While homogenous Nuemann boundary conditions don't alter the process of static condensation, the procedure needs to be adapted to for inhomogenous ones.
First, looking at the definition of the nonlinear residual $R_\mathrm{vol}$ in \eqref{eq:nonlinear_residual:1} we see that this can be split as
\begin{align*}
  R_\mathrm{vol} = R_{\mathrm{vol},\Omega_0} + R_{\mathrm{vol}, \Gamma_{N,0}}
\end{align*}
where $R_{\mathrm{vol}, \Omega_0}$ holds all the terms coming from integration over the domain $\Omega_0$ and $R_{\mathrm{vol}, \Gamma_{N,0}}$ holds all the terms coming from integration over the Neumann surfaces.
Next, note that the bubble functions $\hat{\psi}_\mathrm{B}$ for tetrahedral elements as well as their hexahedral counterparts $\hat{\psi}_{\mathrm{B},1}$, $\hat{\psi}_{\mathrm{B},2}$ have compact support in the finite element interior.
However, the respective gradients don't vanish on the finite element boundary.
Consider an arbitrary finite element \(K \in \mathcal T_h\) with $K \cap \Gamma_{N,0} \neq \emptyset$.
The gradient of a bubble function occurs in bilinear-form \(a_{\Gamma,k}\) in \eqref{eq:cons_lin_neumann},
for the argument \(\Delta \vec{u} \in V_h\).
This yields a non-zero contribution to the element-stiffness-matrix, whereas there is no contribution from
\(R_{\mathrm{vol}, \Gamma_{N,0}}\) to the total element residual vector.
Using the decomposition of local degrees of freedom into \emph{exterior}, $\mathrm{E}$ and \emph{interior}, $\mathrm{I}$
it follows that the local block system can be written in the following form
\begin{align*}
  \begin{pmatrix}
    \tensor{K}_\mathrm{EE} + \tensor{K}_{\Gamma,\mathrm{EE}} & \tensor{K}_\mathrm{EI} + \tensor{K}_{\Gamma,\mathrm{EI}} & \tensor{B}_\mathrm{E} \\
    \tensor{K}_\mathrm{IE} & \tensor{K}_\mathrm{II} & \tensor{B}_\mathrm{I} \\
    \tensor{C}_\mathrm{E} & \tensor{C}_\mathrm{I} & \tensor{D}_\mathrm{E}
  \end{pmatrix}
  \begin{pmatrix}
    \Delta \Rvec{u}_\mathrm{E} \\
    \Delta \Rvec{u}_\mathrm{I} \\
    \Delta \Rvec{p}_\mathrm{E}
  \end{pmatrix}
  =
  \begin{pmatrix}
    -\Rvec{R}_{\mathrm{vol}, \mathrm{E}} - R_{\mathrm{vol},\Gamma, \mathrm{E}} \\
    -\Rvec{R}_{\mathrm{vol},\mathrm{I}} \\
    -\Rvec{R}_{\mathrm{inc},\mathrm{E}}
  \end{pmatrix}.
\end{align*}
The interior degrees of freedom can be statically condensed.
On element level this leads to the static condensed system
\begin{align*}
  \underbrace{
    \begin{pmatrix}
      \widetilde{\tensor{K}} & \widetilde{\tensor{B}} \\
      \widetilde{\tensor{C}} & \widetilde{\tensor{D}}
    \end{pmatrix}
  }_{:= \widetilde{\tensor{A}}}
  \begin{pmatrix}
    \Delta \Rvec{u}_\mathrm{E} \\
    \Delta \Rvec{p}_\mathrm{E}
  \end{pmatrix}
  +
  \underbrace{
    \begin{pmatrix}
      \widetilde{\tensor{K}}_{\Gamma} & \widetilde{\tensor{B}}_{\Gamma} \\
      \tensor 0 & \tensor 0
    \end{pmatrix}
  }_{:= \widetilde{\tensor{A}}_{\Gamma}}
  \begin{pmatrix}
    \Delta \Rvec{u}_\mathrm{E} \\
    \Delta \Rvec{p}_\mathrm{E}
  \end{pmatrix}
  =
  \underbrace{
    \begin{pmatrix}
      -\widetilde{\Rvec{R}}_{\mathrm{vol}}  \\
      -\widetilde{\Rvec{R}}_{\mathrm{inc}}
    \end{pmatrix}
  }_{:= -\Rvec{R}}
  +
  \underbrace{
    \begin{pmatrix}
      - \widetilde{\Rvec{R}}_{\Gamma,\mathrm{upper}} \\
      \Rvec{0}
    \end{pmatrix}
  }_{:= -\Rvec{R}_{\Gamma}},
\end{align*}
where%
\begin{align*}
\widetilde{\tensor{K}}
  &:= \tensor{K}_\mathrm{EE} - \tensor{K}_\mathrm{EI} \tensor{K}_\mathrm{II}^{-1} \tensor{K}_\mathrm{IE},
&\widetilde{\tensor{B}}
  &:= \tensor{B}_\mathrm{E} - \tensor{K}_\mathrm{EI} \tensor{K}_\mathrm{II}^{-1} \tensor{B}_\mathrm{I}, \\
\widetilde{\tensor{C}}
  &:= \tensor{C}_\mathrm{E} - \tensor{C}_\mathrm{I} \tensor{K}_\mathrm{II}^{-1} \tensor{K}_\mathrm{IE},
&\widetilde{\tensor{D}}
  &:= \tensor{D}_\mathrm{E} - \tensor{C}_\mathrm{I} \tensor{K}_\mathrm{II}^{-1} \tensor{B}_\mathrm{I},\\
\widetilde{\tensor{K}}_{\Gamma}
  &:= \tensor{K}_{\Gamma,\mathrm{EE}}
    - \tensor{K}_{\Gamma,\mathrm{EI}} \tensor{K}_\mathrm{II}^{-1} \tensor{K}_\mathrm{IE},
&\widetilde{\tensor{B}}_{\Gamma}
&:= -\tensor{K}_{\Gamma,\mathrm{EI}} \tensor{K}_\mathrm{II}^{-1} \tensor{B}_\mathrm{I},\\
\widetilde{R}_\mathrm{vol}
  &:= \Rvec{R}_{\mathrm{vol},\mathrm{E}}
    - \tensor{K}_\mathrm{EI}\tensor{K}_\mathrm{II}^{-1} \Rvec{R}_{\mathrm{vol},\mathrm{I}},
&\widetilde{R}_{\Gamma,\mathrm{vol}}
  &:= \Rvec{R}_{\Gamma,\mathrm{E}}
    - \tensor{K}_{\Gamma,\mathrm{EI}} \tensor{K}_\mathrm{II}^{-1} \Rvec{R}_{\mathrm{vol},\mathrm{I}},\\
\widetilde{R}_\mathrm{inc}
  &:= \Rvec{R}_{\mathrm{inc},\mathrm{E}}
    - \tensor{C}_\mathrm{I} \tensor{K}_\mathrm{II}^{-1} \Rvec{R}_{\mathrm{inc}, \mathrm{I}}.
\end{align*}
The individual element matrices/vectors $\widetilde{\tensor{A}}$,
$\widetilde{\tensor{A}}_{\Gamma}$,
$\widetilde{\Rvec{R}}$, and $\widetilde{\Rvec{R}}_{\Gamma}$
can be assembled into a global stiffness matrix through loops over volume
elements and surface elements respectively.
In the case of an attached circulatory system a static condensation can be performed in an analogous way
\begin{align*}
  \begin{pmatrix}
    \tensor{K}_\mathrm{EE} + \tensor{K}_{\Gamma,\mathrm{EE}} & \tensor{K}_\mathrm{EI} + \tensor{K}_{\Gamma,\mathrm{EI}} & \tensor{B}_\mathrm{E} & \tensor{E}_{\mathrm{CAV}, \mathrm{E}}\\
    \tensor{K}_\mathrm{IE} & \tensor{K}_\mathrm{II} & \tensor{B}_\mathrm{I} & \tensor 0 \\
    \tensor{C}_\mathrm{E} & \tensor{C}_\mathrm{I} & \tensor{D}_\mathrm{E} & \tensor 0 \\
    \tensor{F}_{\mathrm{CAV}, \mathrm{E}} & \tensor{F}_{\mathrm{CAV}, \mathrm{I}} & \tensor 0 & \tensor{G}_{\mathrm{CAV}}
  \end{pmatrix}
  \begin{pmatrix}
    \Delta \Rvec{u}_{\mathrm E} \\
    \Delta \Rvec{u}_{\mathrm I} \\
    \Delta \Rvec{p}_{\mathrm E} \\
    \Delta \Rvec{p}_{\mathrm{CAV}}
  \end{pmatrix}
  =
  \begin{pmatrix}
    -\Rvec{R}_{\mathrm{vol},\mathrm E} - R_{\Gamma,\mathrm{vol}, \mathrm E} \\
    -\Rvec{R}_{\mathrm{vol}, \mathrm I} \\
    -\Rvec{R}_{\mathrm{inc}, \mathrm E} \\
    -\Rvec{R}_{\mathrm{CAV}, \mathrm E}
  \end{pmatrix}.
\end{align*}
Static condensation of all interior degrees of freedom leads to
\begin{align*}
  \underbrace{
    \begin{pmatrix}
      \widetilde{\tensor{K}} & \widetilde{\tensor{B}} & \tensor{E}_{\mathrm{CAV}, \mathrm E}\\
      \widetilde{\tensor{C}} & \widetilde{\tensor{D}} & \tensor 0 \\
      \widetilde{\tensor{F}}_{\mathrm{CAV}} & \widetilde{\tensor{H}}_{\mathrm{CAV}} & \tensor{G}_{\mathrm{CAV}}
    \end{pmatrix}
  }_{:= \widetilde{\tensor{A}}}
  \begin{pmatrix}
    \Delta \Rvec{u}_{\mathrm E} \\
    \Delta \Rvec{p}_{\mathrm E} \\
    \Delta \Rvec{p}_{\mathrm{CAV}}
  \end{pmatrix}
  +
  \underbrace{
    \begin{pmatrix}
      \widetilde{\tensor{K}}_{\Gamma} & \widetilde{\tensor{B}}_{\Gamma} & \tensor 0\\
      \tensor 0 & \tensor 0 & \tensor 0 \\
      \tensor 0 & \tensor 0 & \tensor 0
    \end{pmatrix}
  }_{:= \widetilde{\tensor{A}}_{\Gamma}}
  \begin{pmatrix}
    \Delta \Rvec{u}_{\mathrm E} \\
    \Delta \Rvec{p}_{\mathrm E} \\
    \Delta \Rvec{p}_{\mathrm{CAV}}
  \end{pmatrix}
  =
  \underbrace{
    \begin{pmatrix}
      -\widetilde{\Rvec{R}}_{\mathrm{vol}}  \\
      -\widetilde{\Rvec{R}}_{\mathrm{inc}}  \\
       \widetilde{\Rvec{R}}_{\mathrm{CAV}}
    \end{pmatrix}
  }_{:= -\Rvec{R}}
  +
  \underbrace{
    \begin{pmatrix}
      - \widetilde{\Rvec{R}}_{\Gamma,\mathrm{vol}} \\
      \Rvec{0} \\
      \Rvec{0}
    \end{pmatrix}
  }_{:= -\Rvec{R}_{\Gamma}},
\end{align*}
where%
\begin{align*}
  \widetilde{\tensor{F}}_{\mathrm{CAV}} &:= \tensor{F}_{\mathrm{CAV}, \mathrm E} - \tensor{F}_{\mathrm{CAV}, \mathrm I} \tensor{K}_\mathrm{II}^{-1} \tensor{K}_\mathrm{IE},\\
  \widetilde{\tensor{H}}_{\mathrm{CAV}} &:= - \tensor{F}_{\mathrm{CAV},\mathrm I}\tensor{K}_\mathrm{II}^{-1}\tensor{B}_{\mathrm I},\\
  \widetilde{\Rvec{R}}_{\mathrm{CAV}} &:= \Rvec{R}_{\mathrm{CAV}} - \tensor{F}_{\mathrm{CAV}, \mathrm I} \tensor{K}_\mathrm{II}^{-1} \Rvec{R}_{\mathrm{vol},\mathrm I}
\end{align*}
%
%
\section{Tensor calculus}
We use the following results from tensor calculus,
for more details we refer to,
e.g.,~\cite{holzapfel2000nonlinear,wriggers2008nonlinear}.
\begin{align*}
  \frac{\partial \overline{\tensor{C}}}{\partial \tensor{C}} &= J^{-\frac{2}{3}} \bb{P} = J^{-\frac{2}{3}} \left(\bb{I} - \frac{1}{3}\tensor{C}^{-1} \otimes \tensor{C}\right),\\
\frac{\partial \tensor{C}^{-1}}{\partial \tensor{C}} &= - \tensor{C}^{-1} \odot \tensor{C}^{-1},\\
{(\tensor{A}\odot\tensor{A})}_{ijkl}&:= \frac{1}{2}\left(A_{ik}A_{jl}+A_{il}A_{jk}\right).
\end{align*}
For symmetric $\tensor{A}$ it holds
\begin{align*}
  \bb{P} : \tensor{A} = \mathrm{Dev}(\tensor{A}) = \tensor{A} - \frac{1}{3}(\tensor{A} : \tensor{C}) \tensor{C}^{-1}.
\end{align*}
%
The isochoric part of the second Piola--Kirchhoff stress tensor as well as the
isochoric part of the fourth order elasticity tensor are given as
\begin{align}
  \label{eq:def_Sisc}\tensor{S}_\mathrm{isc}
    &:= 2 \frac{\partial \overline{\Psi}(\overline{\tensor{C}})}{\partial \tensor{C}}
    = J^{-\frac{2}{3}}\mathrm{Dev}(\overline{\tensor{S}}),\\
  \nonumber\overline{\tensor{S}}
    &:= 2 \frac{\partial \overline{\Psi}(\overline{\tensor{C}})}{\partial \overline{\tensor{C}}},\\
  \label{eq:def_c_isc}\bb{C}_\mathrm{isc}
    & := 4 \frac{\overline{\Psi}(\overline{\tensor{C}})}
                {\partial \tensor{C} \partial \tensor{C}}
      = J^{-\frac{4}{3}} \bb{P} \overline{\bb{C}} \bb{P}^\top
      + J^{-\frac{2}{3}} \frac{2}{3}\mathrm{tr}(\tensor{C} \overline{\tensor{S}})
      \widetilde{\bb{P}}
      - \frac{4}{3}\tensor{S}_\mathrm{isc} \symotimes \tensor{C}^{-1},\\
  \nonumber\overline{\bb{C}}
    &:= 4\frac{\partial \overline{\Psi}
      (\overline{\tensor{C}})}{\partial\overline{\tensor{C}}
      \partial\overline{\tensor{C}}},\\
  \nonumber\widetilde{\bb{P}}
    &:= \tensor{C}^{-1} \odot \tensor{C}^{-1}
      - \frac{1}{3} \tensor{C}^{-1} \otimes \tensor{C}^{-1},\\
  \nonumber\tensor{A}\symotimes \tensor{B}
    &:= \frac{1}{2}\left(\tensor{A}\otimes \tensor{B}
    + \tensor{B} \otimes \tensor{A} \right).
\end{align}

\end{document}